\documentclass[11pt,reqno]{amsart}
\usepackage[letterpaper,textwidth=6.18in,textheight=8.8in,centering]{geometry}
\usepackage[T1]{fontenc}
\usepackage{amssymb,lmodern,microtype}
\usepackage[hidelinks,pdftitle={Invariant flags and integral matrices with prescribed characteristic polynomial}]{hyperref}
\numberwithin{equation}{section}
\newtheorem{thm}{Theorem}[section]
\newtheorem{lem}[thm]{Lemma}
\newtheorem{prop}[thm]{Proposition}
\newtheorem{cor}[thm]{Corollary}
\newcommand{\Oa}{\mathcal O}\newcommand{\R}{\mathbb R}\newcommand{\C}{\mathbb C}\newcommand{\Q}{\mathbb Q}\newcommand{\Z}{\mathbb Z}
\DeclareMathOperator{\vol}{vol}\DeclareMathOperator{\covol}{covol}\DeclareMathOperator{\GL}{GL}\DeclareMathOperator{\SL}{SL}\DeclareMathOperator{\Hom}{Hom}\DeclareMathOperator{\End}{End}
\newcommand{\dd}{\,d}\newcommand{\Frob}{\mathrm F}
\theoremstyle{definition}
\newtheorem{defn}{Definition}[section]

\theoremstyle{remark}
\newtheorem*{example}{Example}
\newcommand{\stage}[1]{\par\smallskip\noindent\emph{#1.}\enspace}
\begin{document}
\title[Invariant flags and integral matrices]{Invariant flags and integral matrices with prescribed characteristic polynomial}

\author{Caihua Luo}
\address{SCHOOL OF SCIENCE AND ENGINEERING,\\
  THE CHINESE UNIVERSITY OF HONG KONG, SHENZHEN,\\
  GUANGDONG 518172, CHINA.}
\email{luocaihua@cuhk.edu.cn}

\author{Xu Song}
\address{SCHOOL OF MATHEMATICAL SCIENCES, HEBEI NORMAL UNIVERSITY,\\
  SHIJIAZHUANG, HEBEI 050024, CHINA.}
\email{songx@hebtu.edu.cn}

\author{Yilin Zeng}
\address{SCHOOL OF SCIENCE AND ENGINEERING,\\
  THE CHINESE UNIVERSITY OF HONG KONG, SHENZHEN,\\
  GUANGDONG 518172, CHINA.}
\email{223015132@link.cuhk.edu.cn}

\date{}
\begin{abstract}
We give an asymptotic formula for integral matrices with a prescribed characteristic polynomial in an expanding entry box over a number field. For cyclic matrices, normalized logarithmic covolumes of a primary invariant flag become uniform on a permutohedral chamber. Their joint distribution with the rescaled matrix factors for tests invariant under reflections at the real places. Separating irreducible factors of multiplicities \(e_i\) into distinct translates multiplies the leading coefficient by \(\prod_i e_i!\). We also bound every noncyclic rational type.
\end{abstract}
\maketitle
\enlargethispage{5pt}
\tableofcontents

\section{Introduction}
How many integral \(n\times n\) matrices in a large box are nilpotent?
Compare this count with the count in the same box of matrices whose eigenvalues are \(0,1,\ldots,n-1\). We shall see that the ratio tends to \(1/n!\). A regular nilpotent matrix has a unique full invariant flag. A matrix with the indicated distinct eigenvalues has one such flag for each ordering of its eigenlines. To pass from this observation to a count of integral matrices one must also determine how the lattices in those flags are distributed.

\subsection{The counting problem}
Let \(K\) be a number field with ring of integers \(\Oa=\Oa_K\). Write
\(S_\infty=\{v:v\mid\infty\}\), \(K_\infty=\prod_v K_v\) (over \(S_\infty\)), \(q_v=[K_v:\R]\), \(q=[K:\Q]=\sum_v q_v\), and \(D_m=\binom m2\), \(D=D_n\).

For \(X=(X_v)_v\in M_n(K_\infty)\), use ordinary real or complex modulus and put
\[
 H_\Box(X)=\max_{v,a,b}|(X_v)_{ab}|,\qquad
 H_{\Frob}(X)=\max_v\|X_v\|_{\Frob},\qquad
 \|X_v\|_{\Frob}^2=\sum_{a,b}|(X_v)_{ab}|^2.
\]
We identify \(M_n(K)\) with its diagonal image. For monic \(f\in\Oa[t]\) of degree \(n\ge1\), let
\[
 \mathcal X_f(\Oa)=\{X\in M_n(\Oa):\chi_X=f\},\qquad
 \mathcal X_f^{\rm cyc}(\Oa)=\{X\in\mathcal X_f(\Oa):\min_X=f\},
\]
where \(\chi_X=\det(tI_n-X)\) and \(\min_X\) is the minimal polynomial over \(K\). The counting function for \(T\ge0\) is
\begin{equation}\label{count}
 N_f^\Box(T)=\#\{X\in\mathcal X_f(\Oa):H_\Box(X)\le T\}.
\end{equation}

Factor \(f=\prod_{i=1}^k p_i^{e_i}\) over \(K\), with distinct monic irreducibles \(p_i\), and put
\[
 d_i=\deg p_i,\qquad n=\sum_i d_i e_i,\qquad
 \ell=\sum_i e_i,\qquad S(T)=T^{qD}(\log T)^{\ell-1}\quad(T>1).
\]
We call \(g\) a \emph{separation of \(f\)} if it is a product
\begin{equation}\label{separate}
 g(t)=\prod_{i=1}^k\prod_{j=1}^{e_i}p_i(t-a_{ij}),\qquad a_{ij}\in\Z ,
\end{equation}
whose displayed factors are pairwise distinct. The shifts are fixed as \(T\to\infty\).
Such shifts exist: the coefficient of \(t^{d_i-1}\) in \(p_i(t-a)\) is its original coefficient minus \(d_i a\). Each previously chosen factor therefore excludes at most one integer \(a\). The \(p_i\) belong to \(\Oa[t]\) since \(\Oa\) is integrally closed.
Our first count can be read directly in terms of \(f\) and \(g\). Theorem~\ref{main} gives
\[
 N_f^\Box(T)\sim \frac{N_g^\Box(T)}{\prod_{i=1}^k e_i!}.
\]
To see where the factorials enter, we now describe a statistic of a cyclic matrix whose distribution distinguishes the two fibers.

\subsection{The invariant flag and its limiting chamber}\label{introflag}
Choose an ordering \(i(1),\ldots,i(\ell)\) in which the \(e_i\) copies of each primary factor occur consecutively. After \(j\) factors the cumulative polynomial and its degree are
\[
 g_j=\prod_{a=1}^j p_{i(a)},\qquad b_j=d_{i(j)},\qquad
 r_j=\deg g_j=\sum_{a=1}^j b_a,\qquad g_0=1,\quad r_0=0.
\]
For a cyclic \(X\), the kernels \(F_j=\ker_K g_j(X)\) form a flag with \(\dim_K F_j=r_j\). Write \(L_j=F_j\cap\Oa^n\) (so \(L_0=0\)) and \(H_j=\covol(L_j)\) in the Euclidean subspace metric of Minkowski space, with \(H_0=1\). The successive covolumes are \(H_j/H_{j-1}\) in the quotient metric. For \(Q(r)=qr(n-r)/2\) and \(T>1\), put
\[
 x_{T,j}=\frac{\log(H_j/H_{j-1})-\ell^{-1}\log H_\ell}{\log T}
                 -Q(r_j)+Q(r_{j-1}).
\]
The logarithms and the \(Q\)-terms telescope, so \(x_T\) belongs to the hyperplane
\(\mathcal H=\{x\in\R^\ell:\sum_jx_j=0\}\).
For \(I\subseteq\{1,\ldots,\ell\}\) write \(b(I)=\sum_{j\in I}b_j\) and \(x(I)=\sum_{j\in I}x_j\). Define
\[
 \begin{aligned}
 Z&=\{x\in\mathcal H:x(I)\ge-Q(b(I))\text{ for all }I\},\\
 P&=\{x\in Z:x_j\le x_{j+1}\text{ if }i(j)=i(j+1)\}.
 \end{aligned}
\]
Integrate on \(\mathcal H\) against \(dx=dx_1\cdots dx_{\ell-1}\), assigning mass one to the point when \(\ell=1\).
The polytope \(Z\) is the centered squarefree polytope of Cai and Deng \cite[Lemma 6.1 and Proposition 5.8]{CD} for the separated degrees. Integral image indices will account for the additional inequalities defining \(P\).

\subsection{The joint distribution and the counting constant}
Let \(\mathcal K_v=O(n)\) for \(K_v=\R\) and \(\mathcal K_v=U(n)\) for \(K_v=\C\). Set \(\mathfrak n_n(K_v)=\{U:U_{ab}=0\text{ for }a\ge b\}\) and
\[
 \eta_{n,v}(B)=\int_{\mathcal K_v}\int_{\mathfrak n_n(K_v)}
           B(kUk^{-1})\,dU\,dk,\qquad \eta_n=\bigotimes_v\eta_{n,v}.
\]
Here \(dk\) is probability Haar measure, and \(dU\) is entry Lebesgue measure, \(dx\,dy\) per complex entry. Write \(r_1=\#\{v:K_v=\R\}\) and \({\rm sgn}:\Oa^\times\to\{\pm1\}^{r_1}\) for the vector of real signs. At a real place let \(R_v\) act as \(\operatorname{diag}(-1,1,\ldots,1)\) there and as the identity elsewhere.
Write \(\kappa_m=\pi^{m/2}/\Gamma(1+m/2)\) for the Euclidean unit ball volume. Unitary or orthogonal invariance of the Frobenius norm gives \(\eta_{n,v}(\mathbf1_{\|\cdot\|_{\Frob}\le1})=\kappa_{q_vD}\).
For a separation \(g\), write \(C^\Frob_{K,g}\) for the squarefree leading coefficient with height \(H_{\Frob}\) in \cite[Theorem 1.2]{CD} and put
\[
 a_{K,g}=C^\Frob_{K,g}\Big/\Big(\vol(Z)\prod_{v\mid\infty}\kappa_{q_vD}\Big).
\]
\begin{thm}[Joint distribution at the free lattice]\label{law}
For fixed \(K,f\), a primary ordering and a separation \(g\), let \(\psi:\mathcal H\to\R\) be bounded continuous and let \(B\in C_c(M_n(K_\infty))\). If \(n\) is even and \({\rm sgn}(\Oa^\times)\ne\{\pm1\}^{r_1}\), assume also
\[
 B(R_vYR_v^{-1})=B(Y)\qquad(v\ \text{real},\ Y\in M_n(K_\infty)).
\]
Then
\begin{equation}\label{lawform}
 \lim_{T\to\infty}\frac{1}{S(T)}
 \sum_{\substack{X\in M_n(\Oa)\\ \min_X=\chi_X=f}}
 \psi(x_T(X))B(X/T)
   = a_{K,g}\int_P\psi(x)\,dx\;\eta_n(B).
\end{equation}

\end{thm}
In particular this allows arbitrary continuous compactly supported angular tests over \(\Q\) and over totally imaginary fields. The sign restriction can be essential, as the quadratic example following Lemma~\ref{rigid} shows. For \(t^2+1\) over \(\Q(\sqrt3)\), the example following Proposition~\ref{localized} computes the mixed-sign proportion \(3/10\). For a cyclic matrix chosen uniformly in the entry box, \(x_T\) tends in distribution to uniform measure on \(P\). The product in \eqref{lawform} expresses its asymptotic independence of the scaled matrix for those tests. Its flag \emph{directions} are a different statistic: at a regular nilpotent limit \(Y_v\), the limiting \(r_j\)-plane is \(\ker Y_v^{r_j}\).

For the entry box we need the local factor \(J_{n,v}=\eta_{n,v}(\mathbf1_{\|\cdot\|_{\max}\le1})\).

\begin{thm}[Counting and separation]\label{main}
For fixed \(K,f\) and every separation \(g\) in \eqref{separate},
\begin{equation}\label{absolute}
 N_f^\Box(T)\sim C_{K,f}^\Box T^{qD}(\log T)^{\ell-1},\qquad
 C_{K,f}^\Box=\frac{C^\Frob_{K,g}}{\prod_i e_i!}
                \prod_{v\mid\infty}\frac{J_{n,v}}{\kappa_{q_vD}} .
\end{equation}
For the single global Euclidean ball the corresponding coefficient is given in Section~\ref{arithsec}. For the same entry box on both fibers,
\begin{equation}\label{relative}
 \lim_{T\to\infty}N_f^\Box(T)/N_g^\Box(T)=1/\prod_i e_i!.
\end{equation}
\end{thm}
For \(n=1\) the count is eventually one, and the zero-dimensional ball has mass one. In particular, if \(g_1,g_2\) are two separations of the same \(f\), then \(C^\Frob_{K,g_1}=C^\Frob_{K,g_2}\). The proof obtains the common density and the ratio \eqref{relative} before using the squarefree counting theorem to evaluate the absolute coefficient.

\subsection{Smaller orbits and examples}
\begin{example}
Take \(K=\Q\) and \(f=t^2\). A nonzero matrix has a unique eigenline, with primitive generator \(v\) up to sign. Then
\[
 H_1=\|v\|_2,\qquad x_{T,1}=\frac{\log\|v\|_2}{\log T}-\frac12,\qquad x_{T,2}=-x_{T,1}.
\]
The full interval for \(x_1\) is \([-1/2,1/2]\). The chamber \(x_1\le x_2\) cuts this to \([-1/2,0]\), and Theorem~\ref{law} gives the uniform law there. At finite height the endpoint can be exceeded. Take \(m\ge2\), \(v=(m,m-1)^t\), \(X=v(1-m,m)\), \(T=m^2\). Then \(X^2=0\), \(H_\Box(X)=T\), and
\[
 x_{T,1}=\frac{\log(1+(1-1/m)^2)}{2\log T}>0.
\]
This gives deviations of order \(1/\log T\) even on regular compact support. Section~\ref{quadratic} computes the limiting distribution directly.
\end{example}

\begin{cor}[The noncyclic contribution]\label{small}
Suppose some $e_i\ge2$ and put $\delta=\min_{e_i\ge2}d_i$. The matrices in \eqref{count} whose minimal polynomial is a proper divisor of $f$ contribute
\[
  O_{K,f}\!\left(T^{q(D_n-\delta)}(\log(2T))^{\ell-2}\right).
\]
\end{cor}

\begin{cor}[The rational coefficient]\label{integercase}
For $K=\Q$, put $L_1=1$, $L_m=\prod_{h=2}^m\pi^{-h/2}\Gamma(h/2)\zeta(h)$, and let $A_p$ be the residue at $s=1$ of the ideal zeta function of $\Z[t]/p$. Then the leading coefficient in \eqref{absolute} is
\[
 \frac{J_{n,\R}}{L_n}\,n^{\ell-2}
       \prod_{i=1}^k\frac{(d_i A_{p_i})^{e_i}}{e_i!}.
\]
\end{cor}

More generally, let \(r_{ij}\) be the number of primary Jordan chains for \(p_i\) of length at least \(j\). The exponent in the upper bound for this rational type is
\[
 \frac q2\bigl(n^2-\textstyle\sum_{i,j}d_i r_{ij}^2\bigr).
\]
The expression in parentheses is the algebraic dimension of its conjugacy orbit. Lemma~\ref{typedimension} computes this dimension. Proposition~\ref{strata} supplies the logarithmic factor as well.

For \(f=(t^2+1)^2\) over \(\Q\), the residue of the ideal zeta function of \(\Z[i]\) is \(\pi/4\) and \(L_4=\pi^2\zeta(3)/1080\). Corollary~\ref{integercase} gives
\[
 N_f^\Box(T)\sim \frac{135J_{4,\R}}{\zeta(3)}\,T^6\log T .
\]
Here the noncyclic type satisfies \(X^2+I=0\) and contributes \(O(T^4)\). Over \(\Q(i)\) the same polynomial has two primary factors, each of multiplicity two. The scale becomes \(T^{12}(\log T)^3\), and the chamber occupies one quarter of that for a separation. Over \(\Q\) it occupies one half. In the split rational cubic case the coefficients for a triple root, a double and a simple root, and three different roots are respectively \(6,18,36\) in units of \(J_{3,\R}/\zeta(3)\).

\subsection{Earlier results and strategy of proof}
Eskin, Mozes and Shah \cite{EMS} related the integral count for an irreducible polynomial over \(\Q\) to its real orbital volume; Shah \cite{Shah} evaluated the arithmetic normalization. For distinct rational roots the logarithm already occurs in the work of Shapira and Zheng \cite{SZ}. Zhang \cite[Theorem 1.5]{Zhang} obtained the squarefree counting asymptotic in a single global Euclidean ball over number fields from translates of maximal tori. His work on observable subgroups \cite{Obs} treats translated homogeneous measures more generally. 

On a repeated primary fiber, even a single rational conjugacy class can contain infinitely many integral classes: the matrices \(mE_{12}\), \(m\ge1\), are separated by their entry ideals. The summation over these classes is part of the fixed-lattice problem. 

Cai and Deng \cite[Theorems 1.1 and 1.2]{CD} compute the squarefree constant in terms of zeta functions of orders; Jeon and Lee \cite{JL} give related orbital integral formulas. We use the normalization of \cite{CD} for \(g\), for which the rational centralizer is a torus. In the cyclic \(f\)-fiber the centralizer is the unit group of \(K[t]/(f)\); its unipotent radical has dimension \(\sum_i d_i(e_i-1)\). Integral conjugacy for such regular matrices is studied by Hertling and Larabi \cite{HL}.

Ostafe and Shparlinski \cite{OS} bound the number of integer matrices with prescribed discriminant. At discriminant zero their count ranges over characteristic polynomials with a repeated root. Here each polynomial, including its multiplicities, is fixed. For \(f\) varying with \(T\), Browning, Gundlach and Sawin \cite{BGS} bound the integer count uniformly in its coefficients. Habegger \cite{Hab} proves an absolute \(O(T^3(\log T)^6)\) bound in degree three. Sarkar \cite{Sarkar} obtains an effective asymptotic in the totally real irreducible cubic case. The constants and limiting operations here concern a fixed fiber. 

For distinct eigenvalues \(\lambda,\mu\), writing \(s=\lambda+\mu\) gives the affine quadric \((2a-s)^2+4bc=(\lambda-\mu)^2\), with \(2a-s\equiv-s\bmod 2\Oa\). Comparisons with Xu--Zhang \cite{XZ}, Kelmer--Kontorovich \cite{KK} and Zhang \cite{quadZ} must therefore retain the congruence and the archimedean region. The double-root fiber is the level-zero cone. Section~\ref{quadratic} treats smooth and double-root fibers by counting eigenlines in the same sharp box.

Here is how the two counting assertions fit together. Each cyclic matrix determines a flag once the primary order is prescribed. Swapping two distinct factors gives two complementary kernels in the same two-layer quotient. Their volume angle is controlled by the resultant. Equal factors instead supply an integral image index and order two consecutive logarithmic heights. These inequalities select \(P\) from \(Z\). Their volume ratio counts integral matrices only if the density above the two chambers agrees.

Section~\ref{entries} computes that density on the space of lattices. Shears fill every rectangular entry except an intertwiner between consecutive equal factors. Its norm tends to infinity in the chamber interior. Averaging its unit orbit supplies ordinary entry measure, uniformly in the remaining shapes. Section~\ref{fixed} bounds the mass of lattices with a nonzero vector near a limiting kernel. This excludes the proper homogeneous measures compatible with the unipotent action and permits localization at \(\Oa^n\). Keeping the prefix-height weight through this argument gives \eqref{lawform}. The singular types are estimated separately in Section~\ref{types}; Section~\ref{arithsec} then identifies the coefficient. Two appendices evaluate the archimedean factors and give an identity at finite height.

\section{Covolumes of invariant subspaces}\label{flags}
Assume \(n\ge2\) throughout the geometric argument. The two flag identities used in Section~\ref{introflag} follow directly from cyclicity and the quotient metric:
\[
 K^n\simeq K[t]/(f),\qquad
 F_j\simeq(f/g_j)/(f),\qquad
 \dim_K F_j=r_j,\qquad
 H_j=H_{j-1}\covol(L_j/L_{j-1}).
\]
The last identity uses the orthogonal projection of \(L_j\) to \(F_j/F_{j-1}\).

Write \(L_{\rm std}=\Oa^n\) in its Minkowski embedding and put
\[
 G=\{h\in\GL_n(K_\infty):\prod_v|\det h_v|^{q_v}=1\},\quad
 \Gamma=\GL_n(\Oa),\quad \mathcal X=G/\Gamma,\quad L=hL_{\rm std} .
\]
Give \(\mathcal X\) its invariant probability \(\mu\). The rational space of \(L\) is \(V_L=L\otimes_{\Oa}K=hK^n\) and its ring of integral endomorphisms is \(h M_n(\Oa)h^{-1}\).
For a cyclic \(X\) on \(V_L\) with polynomial \(f\), let \(A=K[t]/(f)\), acting by \(t\mapsto X\). A cyclic vector identifies \(A\simeq V_L\). For an ideal \(F\lhd A\) use its image under this identification and define
\[
 F_\infty=F\otimes_\Q\R,\qquad h_F(L,X)=\covol_{F_\infty}(L\cap F_\infty).
\]
Two cyclic identifications differ by multiplication by an element of \(A^\times\), which preserves every ideal. Explicitly the ideals are \((d)/(f)\) for the monic divisors \(d\mid f\). Thus \(h_F\) is intrinsic to the transported invariant subspace.

\subsection{How a change of flag changes covolumes}
Let \(\mathcal N=\{Y:Y_v^n=0\ \forall v\}\) and \(\mathcal N^{\rm reg}=\{Y\in\mathcal N:Y_v^{\,n-1}\ne0\ \forall v\}\). A compact matrix set \(\mathcal C\) is \emph{regular} if
\[
        \mathcal C\cap\mathcal N\ \subseteq\ \mathcal N^{\rm reg}.
\]
A test has \emph{regular compact support} if its support is such a set.
Every bounded sequence \(X/T\) from the fixed fiber has only nilpotent limit points as \(T\to\infty\), since \(\chi_{X/T}(t)=T^{-n}f(Tt)\to t^n\) at every place. We abbreviate \(x_{T,j}\) to \(x_j\).

\begin{lem}[Covolumes and image indices]\label{flaggeometry}
Let \(X\) be cyclic, \(L\) an \(X\)-stable lattice as above, and \(X/T\) in a fixed regular compact set. An invariant ideal \(F\) selects an initial segment from each primary chain. If their union of layer indices is \(I\), then
\begin{equation}\label{cuts}
 \log h_F/\log T=Q(b(I))+\sum_{j\in I}x_j+O(1/\log T).
\end{equation}
If layers \(j,j+1\) belong to the same primary, let \(I_j\) be the index of the image of \(p_{i(j)}(X)\) from the upper quotient lattice in the lower quotient lattice. Then
\begin{equation}\label{index}
 \log I_j=(x_{j+1}-x_j)\log T+O(1).
\end{equation}
The implied constants depend on the compact set and on \(K,f\).
\end{lem}
\begin{proof}
\stage{Interchanging two distinct factors}
Interchange adjacent layers of distinct polynomials \(p,\tilde p\) of degrees \(a,b\).
Let \(M\) be the saturated two-layer subquotient and \(M_p,M_{\tilde p}\) its intersections with the two kernels. An integral B\'ezout identity gives
\[
 up+v\tilde p=R,\qquad R=\operatorname{Res}(p,\tilde p)\ne0,\qquad
 RM\subseteq M_p+M_{\tilde p}\subseteq M.
\]
Consequently the integral gluing index is uniformly bounded:
\[
 d_M=[M:M_p+M_{\tilde p}]
       \le[M:RM]=|N_{K/\Q}R|^{a+b}.
\]

At a limit \(X_v/T\to Y_v\) on this compact set, every consecutive two-layer quotient converges along a subsequence to a regular nilpotent block. Indeed, its endpoints converge to the unique Jordan planes of their respective dimensions. A cyclic vector for the limiting quotient supplies a Krylov basis with uniformly bounded condition number on a neighborhood. Finitely many such neighborhoods suffice.

In the resulting polynomial basis put \(p_T(t)=T^{-a}p(Tt)\) and \(\widetilde p_T(t)=T^{-b}\widetilde p(Tt)\), with coefficients at \(v\). The two kernel bases are \(t^i\widetilde p_T\) for \(0\le i<a\) and \(t^j p_T\) for \(0\le j<b\). Their separate wedges tend to wedges of distinct monomials. The combined determinant is the Sylvester determinant up to sign. Thus their real volume angle \(\sigma_v\) satisfies
\[
 \sigma_v\asymp|\operatorname{Res}(p_T,\widetilde p_T)|_v^{q_v}
 =T^{-q_vab}|\operatorname{Res}(p,\widetilde p)|_v^{q_v},
 \qquad \prod_{v\mid\infty}\sigma_v\asymp T^{-qab}.
\]
Here \(|\cdot|_v\) denotes ordinary real or complex modulus. The comparison constants are uniform on the chosen compact set.

Write \(H,H'\) for the two intermediate covolumes and \(H_-,H_+\) for the endpoint covolumes. If \(\sigma\) is the product of the real volume angles and \(d_M=[M:M_p+M_{\tilde p}]\), the subspace and quotient volume identities give
\[
 \frac{HH'}{H_-H_+}=\frac{d_M}{\sigma}\asymp T^{qab},\qquad
 Q(r+a)+Q(r+b)-Q(r)-Q(r+a+b)=qab ,
\]
where \(r\) is the preceding dimension. Let \(x'\) denote the coordinates after the interchange, whose first position is \(j\). The preceding identities give
\[
 (x'_{j+1}-x_j)\log T
     =\log\frac{H_+H_-}{HH'}+qab\log T=O(1).
\]
The sum of the two coordinates is unchanged, so \(x'_j=x_{j+1}+O(1/\log T)\) as well. Move the selected layers to the front using at most \(\binom{\ell}{2}\) interchanges. The sum of their coordinates then telescopes to \(\log h_F/\log T-Q(b(I))+O(1/\log T)\), proving \eqref{cuts}.
\stage{Two copies of one factor}
Write \(d=b_j=b_{j+1}\), \(c_a=H_a/H_{a-1}\), and let \(\overline p_j\) be the map from the upper layer to the lower layer induced by \(p_{i(j)}(X)\). In any convergent sequence of adapted orthonormal frames on a regular compact set, the scaled two-layer operator tends along a subsequence to a regular nilpotent \(J_{2d}\). Writing \(p(t)=t^d+\sum_{h<d}p_ht^h\), we have at each place
\[
 T^{-d}p(X)=(X/T)^d+\sum_{h<d}p_h T^{h-d}(X/T)^h\longrightarrow J_{2d}^d,
 \quad
 J_{2d}^d:K_v^{2d}/\ker J_{2d}^d\xrightarrow{\sim}\ker J_{2d}^d .
\]
Thus each local determinant on the two \(d\)-dimensional layers has size \(\asymp T^{d^2}\) in ordinary modulus; multiplying its real Jacobians gives
\[
 |\det_\R\overline p_j|\asymp T^{q d^2},\qquad
 I_j=[L_j/L_{j-1}:\overline p_j(L_{j+1}/L_j)]
    =|\det_\R\overline p_j|\,\frac{c_{j+1}}{c_j},
\]
where \(L_j=L\cap(F_j)_\infty\). Since \(Q(r+d)-2Q(r)+Q(r-d)=-qd^2\),
\[
 (x_{j+1}-x_j)\log T
 =\log(c_{j+1}/c_j)+qd^2\log T=\log I_j+O(1).
\]
\end{proof}

\begin{cor}[The chamber of admissible heights]\label{chamber}
Fix \(c>0\), a regular compact set \(\mathcal C\), and a compact \(\Omega_0\subset\operatorname{relint} P\). For \(Y=X/T\in\mathcal C\) and all large \(T\) one has
\[
 \begin{aligned}
 h_F\ge c\quad(0<F<A)&\ \Longrightarrow\
          \operatorname{dist}(x_T,P)\le C/\log T,\\
 x_T\in \Omega_0 &\ \Longrightarrow\
          h_F\ge T^{\gamma}\quad(0<F<A)
 \end{aligned}
\]
for \(C=C(K,f,c,\mathcal C)>0\) and some \(\gamma>0\) depending also on \(\Omega_0\). For a separation the chamber is \(Z\), and
\begin{equation}\label{ch}
 \vol(P)=\vol(Z)/\prod_i e_i!.
\end{equation}
\end{cor}
\begin{proof}
Fix integers \(0\le m_i\le e_i\). Let \(\mathcal A(\boldsymbol m)\) be the subsets selecting \(m_i\) indices from run \(i\), and let \(I_0\) select the first \(m_i\) there. By \eqref{index}, \(I_j\ge1\) gives \(x_j\le x_{j+1}+O(1/\log T)\) inside a run. Consecutive exchanges therefore give
\[
 b(I)=b(I_0),\qquad x(I)\ge x(I_0)-O(1/\log T)
          \quad(I\in\mathcal A(\boldsymbol m)).
\]
The ideal for \(I_0\) satisfies \eqref{cuts}. Its lower height bound gives \(x(I_0)+Q(b(I_0))\ge -O(1/\log T)\), which now holds for all \(I\) in the family. The fixed finite system of chamber inequalities gives the asserted distance bound by the linear error bound for a polyhedron (or by triangulating its normal cones). Conversely, on a compact set in the interior, \(x(I)+Q(b(I))\) has a positive minimum over all nonempty proper \(I\). Taking \(\gamma\) below that minimum and using \eqref{cuts} gives \(h_F\ge T^\gamma\). When \(\ell=1\) there are no proper ideals.

With \(\mathfrak S=\prod_i S_{e_i}\) acting on the coordinates of \(Z\), let \(P^\circ=\operatorname{relint}P\). A small perturbation of \(0\) strictly orders the equal-primary coordinates and keeps every proper subset inequality strict. For \(\ell>1\) this proves full relative dimension, and
\[
  Z=\bigcup_{\sigma\in\mathfrak S}\sigma P,\qquad
  \sigma P^\circ\cap\tau P^\circ=\varnothing\quad(\sigma\ne\tau),\qquad
  \vol Z=|\mathfrak S|\vol P .
\]
All these permutations preserve \(dx\), and the boundaries have zero \(dx\)-measure. For \(\ell=1\) the identity holds with the point-mass convention.
\end{proof}

\subsection{Unfolding the flag}
For a flag of \(K\)-subspaces \(F_\bullet\), put \(L_j=L\cap(F_j)_\infty\) and \(P_j=L_j/L_{j-1}\).

\begin{lem}[Integral flags and their measure]\label{flagmeasure}
Fix positive layer dimensions \(b_1,\ldots,b_\ell\) with sum \(n\). The \(\Gamma\)-orbits of flags on \(L_{\rm std}\) are parametrized by the ordered tuples of projective \(\Oa\)-modules satisfying
\[
 \operatorname{rank}_\Oa P_j=b_j,\qquad \bigoplus_jP_j\simeq\Oa^n .
\]
There are finitely many such tuples. For each tuple the unfolded flag measure, in coordinates \(k\operatorname{diag}(a_jS_j)u\), has height factor
\begin{equation}\label{jac}
 c_{P_\bullet}\prod_{j<l}(a_j/a_l)^{qb_jb_l}\prod_{j<\ell}d\log H_j .
\end{equation}
Here \(a_j>0\), \(\prod_j a_j^{b_j}=1\), \(|\det_\R S_j|=1\), and \(c_{P_\bullet}>0\). The other factors are compact frame probability, the independent normalized block-shape measures, and probability on the shear quotient.
\end{lem}
\begin{proof}
Saturation and the Dedekind property give
\[
 L_j/L_{j-1}\text{ torsion-free}\ \Longrightarrow\
 P_j\text{ projective}\ \Longrightarrow\
 L_j\simeq L_{j-1}\oplus P_j .
\]
Iterating gives a splitting compatible with the full flag. Isomorphisms of corresponding \(P_j\) extend by direct sum to an automorphism of \(L_{\rm std}\), which proves the classification. By the Steinitz classification,
\[
 P_j\simeq\Oa^{b_j-1}\oplus\mathfrak a_j,\qquad
 \bigoplus_jP_j\simeq\Oa^n
 \ \Longleftrightarrow\ \prod_j[\mathfrak a_j]=1
       \quad\text{in }\operatorname{Cl}(\Oa).
\]
Finiteness follows from \(|\operatorname{Cl}(\Oa)|<\infty\).

Choose the fractional ideal representatives so that \(\prod_j\mathfrak a_j=\Oa\), by rescaling one of them by the inverse of a generator of their principal product.
For the fixed block model \(P_j=\Oa^{b_j-1}\oplus\mathfrak a_j\), let \(h_j^0\) be its Minkowski covolume. Then
\[
 \prod_j h_j^0
   =\covol(\Oa)^n\prod_j N\mathfrak a_j
   =\covol(L_{\rm std})=H_\ell .
\]
Thus the fixed direct-sum model has the ambient covolume. Since each \(S_j\) has real determinant of modulus one, its positive block scales satisfy
\[
 \prod_j a_j^{qb_j}=1\quad\Longleftrightarrow\quad
 \prod_j a_j^{b_j}=1.
\]

In the \(jl\) rectangle conjugation by \(a=\operatorname{diag}(a_j)\) has real determinant \((a_j/a_l)^{qb_jb_l}\). Thus the parabolic modulus in the order \(kaSu\) is their product. The prefix covolumes satisfy
\[
 H_j=\prod_{h\le j}h_h^0 a_h^{qb_h},\qquad
 d\log H_j=\sum_{h\le j}qb_h\,d\log a_h .
\]
On \(\sum_hb_h\log a_h=0\), this linear change of variables has a fixed nonzero Jacobian in the first \(\ell-1\) coordinates, absorbed in \(c_{P_\bullet}\). Each norm-one block quotient has finite invariant mass by reduction theory and the unit theorem. The single global determinant constraint is exhausted by \(\prod a_j^{b_j}=1\), leaving the block shapes independent. Normalizing the finite masses gives \eqref{jac}.
\end{proof}

\section{Equidistribution above the chamber}\label{entries}
\subsection{Integration over the shears}
Fix a tuple \(P_\bullet\) in Lemma~\ref{flagmeasure}. The diagonal models are endomorphisms \(A_j\in\End_\Oa(P_j)\) with \(\chi_{A_j}=p_{i(j)}\). A field frame identifies \(P_j\) with a full lattice in \(E_j=K[t]/p_{i(j)}\). Write
\[
 \begin{aligned}
 E_{j,\infty}^1&=\{z\in E_{j,\infty}^\times:|N_{E_j/\Q}z|=1\},&
 U_j&=\{z\in E_j^\times:zP_j=P_j\},\\
 \Lambda_{jl}&=\Hom_\Oa(P_l,P_j),&
 M_{jl}&=\Lambda_{jl}\otimes_\Z\R .
 \end{aligned}
\]
The shear group \(\mathcal U\) consists of block upper unitriangular matrices with \(u_{jl}\in M_{jl}\). Its integral subgroup \(\mathcal U_\Z\) has \(u_{jl}\in\Lambda_{jl}\); give \(\mathcal U/\mathcal U_\Z\) probability Haar measure.

\begin{lem}[Finite diagonal models and translation]\label{diagonalmodels}
The diagonal tuples \(A_\bullet\) have finitely many classes under \(\prod_j\operatorname{Aut}_\Oa(P_j)\). For each tuple, every \(E_{j,\infty}^1/U_j\) is compact. The correspondence
\[
 A_j\longmapsto A_j+a_{ij}I
\]
between a fixed primary factor and its translate preserves the integral centralizer, its compact-quotient volume, and the entry lattice \(\Lambda_{jl}\).
\end{lem}
\begin{proof}
For \(E=K[t]/p_i\), \(R=\Oa[t]/p_i\), and a full \(R\)-lattice \(I\subset E\), its maximal-order saturation \(J=\Oa_E I\) satisfies
\[
 \mathfrak fJ\subseteq I\subseteq J,\qquad \mathfrak f=(R:\Oa_E).
\]
After scaling, \(J\) ranges over a finite set \(\mathcal J_E\) of ideal-class representatives. Thus
\[
 \#\big(\{I\}/E^\times\big)
 \le\sum_{J\in\mathcal J_E}
      \#\{M\le J/\mathfrak fJ:\ M\text{ an }R\text{-submodule}\}<\infty.
\]
The subcollection whose underlying \(\Oa\)-module is \(P_j\) is finite as well.
The multiplier order gives
\[
 R\subseteq(I:I)\subseteq\Oa_E,\qquad
 U_I=(I:I)^\times,\qquad[\Oa_E^\times:U_I]<\infty.
\]
The unit theorem and this finite index imply compactness of \(E_\infty^1/U_I\).
Finally
\[
 h(A+aI)=(A+aI)h\ \Longleftrightarrow\ hA=Ah,\qquad
 \Oa[A+aI]=\Oa[A].
\]
The underlying projective layers also agree. These identities prove all three invariances.
\end{proof}
In every unfolded model, give each compact centralizer quotient mass one and absorb its original volume into the model multiplier.

Let \(A\) have diagonal blocks \(A_j\) and upper entries \(A_{jl}\in\Lambda_{jl}\); set \(C=uAu^{-1}\). On \(M_{jl}\) write
\[
 \begin{aligned}
 \mathcal S_{jl}v&=vA_l-A_jv,& V_{jl}&=\ker \mathcal S_{jl},& W_{jl}&=\operatorname{im}\mathcal S_{jl},\\
 \pi_{jl}&:M_{jl}=V_{jl}\oplus W_{jl}\longrightarrow V_{jl},&
 I_{jl}&=\pi_{jl}\Lambda_{jl},& K_{jl}&=V_{jl}\cap\Lambda_{jl}.
 \end{aligned}
\]
The projection is along \(W_{jl}\), and measures satisfy \(dC=dz\,dw\) in \(V_{jl}\oplus W_{jl}\). The exceptional edges are \(\mathcal E=\{(j,j+1):i(j)=i(j+1)\}\).
Semisimplicity gives the displayed direct sum: after extending scalars, \(\mathcal S_{jl}\) acts on each spectral block by an eigenvalue difference. Thus \(V_{jl}=0\) for distinct primaries, and \(V_{jl}\simeq E_{j,\infty}\) in fixed field frames for equal primaries. The two lattices \(I_{jl}\) and \(K_{jl}\) will play different roles in the proof.

\begin{lem}[Distribution of the upper blocks]\label{sheardist}
Fix the scales and the diagonal blocks \(A_j\). Sum over cyclic block upper triangular matrices \(A\) with these diagonal blocks and \(A_{jl}\in\Lambda_{jl}\) for \(j<l\). For every nonnegative Borel function \(\Phi\) of the upper entries,
\[
 \int_{\mathcal U/\mathcal U_{\Z}}\ \sum_{A\ {\rm cyclic}}\Phi((uAu^{-1})_{jl})\,du
   =\int\Phi((C_{jl}))\prod_{j<l}d\nu_{jl}(C_{jl}).
\]
For a single block abbreviate \((\Lambda_{jl},W_{jl},I_{jl})\) to \((\Lambda,W,I)\). Its measure is determined for nonnegative \(\phi\) by
\begin{equation}\label{exc}
 \nu_{jl}(\phi)=
 \begin{cases}
 \displaystyle\frac{1}{\covol\Lambda}\int\phi(C)\,dC,
                      &(j,l)\notin\mathcal E,\\[6pt]
 \displaystyle\frac{1}{\covol(\Lambda\cap W)}
       \sum_{z\in I\setminus\{0\}}\int_W\phi(z+w)\,dw,
                      &(j,l)\in\mathcal E .
 \end{cases}
\end{equation}
Heights of the prefix subspaces may be held fixed in this identity.
\end{lem}

\begin{proof}
\stage{One block}
For \(\mathcal S=\mathcal S_{jl}\) one has \(\mathcal S\Lambda\subseteq\Lambda\).
Integration over the torus of shear entries pushes its probability to Haar probability on \(W/(W\cap\Lambda)\). Indeed the image is a compact connected subgroup with Lie algebra \(W\). If \(t\) is the offset from shorter blocks, then for any nonnegative test \(\phi\),
\[
 \int\sum_{b\in\Lambda}\phi(b+\mathcal S v+t)\,dv
 =\frac1{\covol(W\cap\Lambda)}
       \sum_{z\in I}\int_W\phi(z+\pi t+w)\,dw .
\]
Indeed \(\pi^{-1}(z)\cap\Lambda=b_z+(\Lambda\cap W)\) for any lift \(b_z\), and for a period cell \(\mathcal P_W\) one has
\[
 \int_{\mathcal P_W}\sum_{\lambda\in\Lambda\cap W}\phi(b_z+t+w+\lambda)\,
             \frac{dw}{\covol(\Lambda\cap W)}
   =\frac{1}{\covol(\Lambda\cap W)}\int_W\phi(z+\pi t+w)\,dw .
\] To disintegrate the unused shear coordinates, write \(M=M_{jl}\), \(V=V_{jl}\) and \(K_0=V\cap\Lambda\). There is an exact sequence of tori
\[
 0\longrightarrow V/K_0\longrightarrow M/\Lambda
   \longrightarrow (M/V)/\operatorname{pr}_{M/V}\Lambda\longrightarrow0 .
\]
Condition on the last coordinate. Each fiber is a translate of \(V/K_0\) with Haar probability. Conditioning only on the value of \(\mathcal S v\bmod(W\cap\Lambda)\) can leave several such fibers. Indeed, for the induced map \(\overline{\mathcal S}:M/\Lambda\to W/(W\cap\Lambda)\),
\[
 \ker\overline{\mathcal S}\,/\,(V/K_0)\ \simeq\ (W\cap\Lambda)/\mathcal S\Lambda,\qquad
 [v]\longmapsto [\mathcal S v] .
\]
For diagonal entries \(0,2\) over \(\Z\) this finite group has order two. Haar probability on the target still has total mass one. The higher-block calculation below conditions on the transverse coordinate in the exact sequence.

\stage{The first superdiagonal}
Here the offset is zero. For \(p(t)=\sum_k p_k t^k\) on two equal-primary layers with diagonal blocks \(A_j,A_l\), the induced upper block is
\[
 Dp_{A_j,A_l}(A_{jl})=\sum_{k\ge1}p_k\sum_{r=0}^{k-1} A_j^r A_{jl} A_l^{k-1-r}.
\]
For \(A_{jl}=vA_l-A_jv\) this telescopes to \(vp(A_l)-p(A_j)v=0\). For \(A_jA_{jl}=A_{jl}A_l\) it equals \(p'(A_j)A_{jl}\). Thus \(Dp(A_{jl})=p'(A_j)\pi A_{jl}\). On a run of \(e\) copies of \(p\), identify the layers with \(E=K[t]/p\) and put \(\alpha=t\bmod p\). Each application of \(p(A)\) lowers the filtration by at least one step. For \(e\ge2\), the map of the last layer to the first after \(e-1\) steps is therefore multiplication by
\[
 (p'(\alpha))^{e-1}z_1\cdots z_{e-1},
 \qquad z_j=\pi A_{j,j+1}\in E .
\]
Since \(p'(\alpha)\ne0\) and \(E\) is a field, the primary minimal polynomial is \(p^e\) exactly when every \(z_j\ne0\). For \(e=1\) this condition is empty. Distinct primary components split over \(K\) by B\'ezout's identity. Applying the criterion to each run gives \eqref{exc}.

\stage{The higher superdiagonals}
Order the blocks by distance \(h=l-j\). Use additive entry parallelepipeds for \(u\), a fundamental domain by successive triangular reduction. Eliminate the output variables in \emph{decreasing} distance. The identity \(Cu=uA\) gives
\[
 C_{jl}=A_{jl}+u_{jl}A_l-A_ju_{jl}
           +\sum_{j<h<l}(u_{jh}A_{hl}-C_{jh}u_{hl}).
\]
At distance \(h\ge2\), first integrate the distance-\(h\) shear variables as in the one-block calculation. What remains for equal factors is a periodic function of the projected offsets in the torus
\(\mathbb T_h=\prod_j V_{j,j+h}/I_{j,j+h}\), the product running over the edges at this distance within a single run. Keep all shorter integral entries and the transverse parts of their shears fixed. The unused kernel coordinates \(v_j\in V_{j,j+h-1}\) at distance \(h-1\) leave all shorter outputs fixed. By the displayed recursion they translate the offsets by
\[
 (L_hv)_j=v_jz_{j+h-1}-z_jv_{j+1}.
\]
For an integral intertwiner \(v\), composition by \(v\) preserves both the kernel and image summands of the Sylvester operator. For every conformable block \(Z\) this gives \(\pi(vZ)=v\pi(Z)\) and \(\pi(Zv)=\pi(Z)v\). In particular,
\[
 L_h\Big(\prod_j K_{j,j+h-1}\Big)\ \subseteq\ \prod_j I_{j,j+h}.
\]
The induced map of tori is onto. Given any \(y=(y_j)\) in the target vector space, set the last source variable to zero and solve, from right to left,
\[
              v_j=(y_j+z_jv_{j+1})z_{j+h-1}^{-1}.
\]
Every \(z_{j+h-1}\in E^\times\) is invertible at every infinite place. The real linear map is therefore onto and so is its quotient on tori.

It follows that for every nonnegative periodic \(F\), with the shorter data held as specified,
\[
 \int F(t+L_hv)\,dv=\int_{\mathbb T_h}F(\xi)\,d\xi
\]
in probability measures. Unfolding this last torus gives \(dC/\covol\Lambda\), since
\(\covol\Lambda=\covol(W\cap\Lambda)\covol I\).
This eliminates that entire distance jointly. The kernel variables just used have no effect on the next shorter outputs, and the longer outputs have already been integrated. Iteration leaves precisely the adjacent measures above. Tonelli's theorem permits all these operations for nonnegative tests.
\end{proof}

\begin{example}[Three copies of \(t^2+1\)]
In the standard lattice take \(J=\left(\begin{smallmatrix}0&-1\\1&0\end{smallmatrix}\right)\). The projection on commuting matrices is \(Z\mapsto(Z-JZJ)/2\); its integral image is \(\frac12\Z[i]\), whereas the integral commuting matrices form \(\Z[i]\). With three copies the sole higher offset is translated by \(v_{12}z_{23}-z_{12}v_{23}\). The variables \(v_{12},v_{23}\) have periods \(\Z[i]\), and the target has periods \(\frac12\Z[i]\). For nonzero \(z_{23}\), multiplication by \(z_{23}\) maps the first torus onto the target. Its pushforward has total mass one, irrespective of its finite kernel. This is the reason that a third lattice index does not appear.
\end{example}

\subsection{Averaging the entries on an equal-primary edge}
For a number field \(E\) write \(m=[E:\Q]\) and \(E_\infty^1=\{h\in E_\infty^\times:|Nh|=1\}\).
\begin{lem}[Averaging a field lattice]\label{shell}
Let \(I\subset E\) be a full lattice stable under a finite-index subgroup \(U\le\Oa_E^\times\), and let \(\varepsilon>0\). With \(dh\) probability on \(E_\infty^1/U\), define
\[
 \nu_{I,\varepsilon}=
 \frac{1}{\log(1+\varepsilon)}
 \int_1^{1+\varepsilon}\int_{E_\infty^1/U}\sum_{0\ne z\in I}\delta_{shz}\,dh\,\frac{ds}{s},
 \qquad
 b_{m,\varepsilon}=\frac{1-(1+\varepsilon)^{-m}}{m\log(1+\varepsilon)} .
\]
Then \(\nu_{I,\varepsilon}=\rho_I\,dx\), where
\[
 \lim_{M\to\infty}\operatorname*{ess\,sup}_{|Nx|\ge M}
       |\covol(I)\rho_I(x)-b_{m,\varepsilon}|=0,\qquad
 \|\rho_I\|_\infty\le \frac{C}{\covol(I)} .
\]
The density bound is uniform over finitely many homothety classes of \(I\) for fixed \(U,\varepsilon\). The convergence is asserted for each fixed lattice.
\end{lem}
\begin{proof}
Write \(H=E_\infty^1\). Normalize Haar measure \(dh_0\) on \(H\) so that \(H/U\) has mass one for counting measure on \(U\). In coordinates \(x=t^{1/m}h\), \(t=|Nx|>0\), scaling additive volume gives \(dx=\kappa_U\,dt\,dh_0\) for some \(\kappa_U>0\). Let \(A_I(y)=\#\{[z]\in(I\setminus0)/U:|Nz|\le y\}\) and \(\alpha=(1+\varepsilon)^{-m}\). Unfolding each orbit, and putting \(t=s^m|Nz|\), gives
\[
 \frac{ds}{s}=\frac{dt}{mt},\qquad
 \rho_I(x)=\frac{A_I(t)-A_I(\alpha t)}
               {m\kappa_U t\log(1+\varepsilon)}
 \quad\text{for almost every }x .
\]
Indeed the shell from an orbit of norm \(y\) is \(y\le t\le y/\alpha\); its endpoints are null.

Choose a bounded logarithmic unit parallelepiped and representatives for its finite torsion action on the signs and phases. Restricting \(0<t\le1\) to these representatives gives a bounded fundamental region \(\mathcal F\) with null boundary and \(\vol\mathcal F=\kappa_U\). The radii at the infinite places are all bounded multiples of \(t^{1/m}\), so the only extra limiting point at \(t=0\) is the origin. Dilation of this region gives
\[
 A_I(y)=\#\big((I\setminus0)\cap y^{1/m}\mathcal F\big)
        =\frac{\kappa_U}{\covol I}\,y+o_I(y).
\]
Substitution gives \(\covol(I)\rho_I(x)\to(1-\alpha)/(m\log(1+\varepsilon))\) uniformly away from shell endpoints as \(t\to\infty\). Since \(I\) has bounded denominators, the nonzero absolute norms on \(I\) have a positive lower bound. The shell formula is therefore bounded also for bounded \(t\). Finally
\(\rho_{\lambda I}(x)=|N\lambda|^{-1}\rho_I(x/\lambda)\) and
\(\covol(\lambda I)=|N\lambda|\covol I\) for \(\lambda\in E^\times\), which proves the asserted uniform bound.
\end{proof}

\subsection{The continuous integral and comparison of edge measures}
First compute the scaled measure on one diagonal block. Then combine the blocks and integrate the flag heights. In both lemmas \(\|\cdot\|\) denotes the product Frobenius norm and \(\eta_b\) denotes the compact average of strictly upper triangular entry measure in size \(b\), as in the introduction.

\begin{lem}[Limit of a simple diagonal orbit]\label{simpleorbital}
Fix a diagonal model \(A\) of size \(b\) with irreducible characteristic polynomial over \(K\). Let \(d\omega_A(S)\) be its unfolded shape measure, including the compact arithmetic centralizer quotient, and let \(dk_b\) be compact frame probability. Put
\[
 \sigma_{A,T}(\phi)=T^{-qD_b}\int
    \phi(k_b S A S^{-1}k_b^{-1}/T)\,dk_b\,d\omega_A(S).
\]
There is \(c_A>0\) such that \(\sigma_{A,T}\to c_A\eta_b\) vaguely. More precisely, for \(T\ge2\) and Lipschitz \(\phi\) supported on \(\{\|Y\|\le M\}\),
\[
 \big|\sigma_{A,T}(\phi)-c_A\eta_b(\phi)\big|
       \ll_{A,M} T^{-1}\big(\|\phi\|_\infty+\operatorname{Lip}\phi\big).
\]
\end{lem}
\begin{proof}
\stage{The spectral rectangles} At a complex place choose one-dimensional spectral blocks. At a real place choose a real line for each real root and a real plane for each conjugate pair. These spectral blocks have disjoint spectra. In local coordinates \(kmu'=kum\), conjugation \(u=mu'm^{-1}\) absorbs the parabolic modulus in \(du\). More explicitly, if \(P=MN\) is the block parabolic and \(\delta_P(m)=|\det(\operatorname{Ad}(m)|_{\operatorname{Lie}N})|\), then
\[
 du=\delta_P(m)\,du',\qquad
 dk\,\delta_P(m)\,dm\,du'=dk\,du\,dm
\]
up to a fixed constant. The centralizer in the global determinant-modulus-one group has compact arithmetic quotient. Positive scalar multiples can impose the global modulus constraint without changing a conjugate. These choices contribute fixed factors to the orbital measure.

For the resulting diagonal spectral blocks \(B_r\) of sizes \(h_r\), the identity \(Cu=u\operatorname{diag}(B_r)\) gives, for \(r<s\),
\[
 C_{rs}=u_{rs}B_s-B_ru_{rs}-\sum_{r<j<s}C_{rj}u_{js}.
\]
In increasing order of \(s-r\) this is triangular with diagonal
\(\mathcal S_{rs}:V\mapsto VB_s-B_rV\).
Over \(K_v\) its determinant is \(\prod_{\lambda\in\mathrm{spec} B_r,\ \mu\in\mathrm{spec} B_s}(\mu-\lambda)\), with algebraic multiplicities. It is nonzero and constant on the block orbits. Thus changing from \(u_{rs}\) to \(C_{rs}/T\) supplies Lebesgue measure and the power \(T^{q_v h_rh_s}\) per rectangle.

\stage{A real conjugate pair} For a real pair \(\alpha\pm i\beta\), put \(a=|\beta|/T\) and \(x=ae^{2t}\), \(t\ge0\). The polar integration formula (see \cite[(3.7)--(3.8)]{GOS}) gives
\[
 T^{-1}\sinh(2t)\,dt=\frac{1}{4|\beta|}\,
       \mathbf1_{x\ge a}(1-a^2/x^2)\,dx .
\]
The diagonal entries of the scaled block are \(\alpha/T\); its off-diagonal entries are \(-\operatorname{sgn}(\beta)x\) and \(\operatorname{sgn}(\beta)a^2/x\). Hence the displacement from its limiting triangular block is \(O_A(T^{-1})\) since \(a^2/x\le a\). The radial density has error
\[
 \int_0^\infty\big|1-\mathbf1_{x\ge a}(1-a^2/x^2)\big|\,dx
      =a+\int_a^\infty a^2x^{-2}\,dx=2a .
\]
An orthogonal reflection supplies the opposite sign under compact averaging; the intervening rectangles have invariant Lebesgue measure. Real singletons and complex diagonal eigenvalues have displacement \(O_A(T^{-1})\) as well.

\stage{The orbital limit} The local exponents add to
\[
 \begin{aligned}
 \sum_{r<s}h_rh_s+c&=D_b
       &&(v\text{ real, }c\text{ conjugate pairs}),\\
 2\sum_{r<s}1&=2D_b &&(v\text{ complex}),\\
 \sum_v q_vD_b&=qD_b .
 \end{aligned}
\]
The limiting coordinates fill a strictly upper triangular matrix in a compact frame.
On a fixed matrix ball all free coordinates are bounded.
Telescoping the radial densities and then moving the test argument gives
\[
 |\sigma_{A,T}(\phi)-c_A\eta_b(\phi)|
 \le \underbrace{C_{A,M}T^{-1}\|\phi\|_\infty}_{\text{density error}}
     +\underbrace{C_{A,M}T^{-1}\operatorname{Lip}\phi}_{\text{displacement error}} .
\]
All Jacobian constants are positive. The same coordinates give
\(\sup_{T\ge2}\sigma_{A,T}(\{\|Y\|\le M\})<\infty\).
Uniform Lipschitz approximation on compact sets now proves vague convergence.
\end{proof}

Fix the tuple of layers and diagonal models, denoted by \(\mathfrak m\). In the flag coordinates of \eqref{jac} write \(aS=\operatorname{diag}(a_jS_j)\), \(C_{jj}=A_j\), \(C_{lj}=0\) for \(j<l\). For a product of entry measures \(d\nu\) put
\[
 \begin{aligned}
 Y_T&=T^{-1} k(aS)C(aS)^{-1}k^{-1},\\
 \mathcal I_T(B;\Omega_1;d\nu)&=\int\mathbf1_{\Omega_1}(x_T)B(Y_T)\,d\nu\,
   \prod_{j<l}(a_j/a_l)^{qb_jb_l}\prod_{j<\ell}d\log H_j\,dk\,\prod_j d\omega_j .
 \end{aligned}
\]
A fixed model multiplier is understood. Use \(\mathcal I_T\) for \(d\nu=\prod d\nu_{jl}\) from \eqref{exc}, and \(\mathcal I_T^0\) for \(d\nu=\prod dC_{jl}/\covol\Lambda_{jl}\). No ideal cuts have yet been imposed.

\begin{lem}[Factorization of the continuous integral]\label{continuouslemma}
For each \(\mathfrak m\) there are locally finite positive measures \(\beta_{T,\mathfrak m}\), independent of the heights, such that for bounded measurable \(\Omega_1\subset\mathcal H\) and bounded Borel \(B\) of compact support,
\begin{equation}\label{continuousfactor}
 \mathcal I_T^0(B;\Omega_1)=S(T)|\Omega_1|\beta_{T,\mathfrak m}(B).
\end{equation}
Moreover, \(\beta_{T,\mathfrak m}\to c_{\mathfrak m}\eta_n\) vaguely for some \(c_{\mathfrak m}>0\), and
\(\sup_{T\ge2}\beta_{T,\mathfrak m}(\{\|Y\|\le M\})<\infty\).
When \(\ell=1\), \(dx\) is point mass on the single height.
\end{lem}
\begin{proof}
For \(j<l\) write \(Z_{jl}=T^{-1}(a_j/a_l)S_jC_{jl}S_l^{-1}\). Since \(|\det_\R S_j|=|\det_\R S_l|=1\), the height-dependent Jacobians cancel:
\[
 \prod_{j<l}(a_j/a_l)^{q b_j b_l}\,dC_{jl}
       =T^{q\sum_{j<l}b_jb_l}\prod_{j<l}dZ_{jl}.
\]
For \(j<\ell\) the centered flag coordinates satisfy
\[
 \log H_j=(\log T)\big(\textstyle\sum_{k\le j}x_k+Q(r_j)\big)
                  +(j/\ell)\log H_\ell,\qquad
 \prod_{j<\ell}d\log H_j=(\log T)^{\ell-1}\,dx.
\]
The derivative matrix in \(x_1,\ldots,x_{\ell-1}\) is lower triangular. Since \(D_n=\sum_{j<l}b_jb_l+\sum_jD_{b_j}\), what remains after division by \(S(T)\) is \(dx\) times the image measure on \(Y=kZk^{-1}\) specified by
\[
 Z_{jj}=T^{-1}S_jA_jS_j^{-1},\qquad
 dk\,\prod_j T^{-qD_{b_j}}d\omega_j(S_j)
          \prod_{j<l}\frac{dZ_{jl}}{\covol\Lambda_{jl}} .
\]
Include the fixed model multiplier in this image measure to obtain \(\beta_{T,\mathfrak m}\) and \eqref{continuousfactor}.

Insert independent compact rotations on the diagonal layers. The simultaneous changes \(Z_{jl}\mapsto k_j Z_{jl} k_l^{-1}\) and \(k\mapsto k\operatorname{diag}(k_j)^{-1}\) preserve the measures. Lemma~\ref{simpleorbital} now applies independently to each diagonal block. On a bounded ball all the rectangles are bounded and each diagonal measure has uniformly bounded mass. Products therefore converge on continuous compact tests (first use finite sums of product tests, then uniform approximation on the bounding compact set). The within-layer and between-layer entries fill all strictly upper triangular coordinates. This identifies the limit as \(c_{\mathfrak m}\eta_n\) and gives the asserted ball bound.
\end{proof}

The identity \eqref{continuousfactor} holds for Borel tests at each \(T\). The angular limit requires a continuity condition. For example, if \(K=\Q\), \(A\) has polynomial \(t^2+1\), and \(\Omega=\{Y:Y^2=0,\ \|Y\|_{\Frob}\le1\}\), then
\[
 \det(SAS^{-1}/T)=T^{-2}\quad\Longrightarrow\quad
 \sigma_{A,T}(\Omega)=0,\qquad \eta_2(\Omega)=\int_{-1}^{1}du=2 .
\]
We now compare the original edge measures with the continuous integral on the chamber interior.

\begin{lem}[Replacing the adjacent edge measures]\label{joint}
For \(\Omega_0\subset\operatorname{relint}P\) compact and \(B\in C_c\) of regular compact support,
\[
  \lim_{T\to\infty}
  \frac{\mathcal I_T(B;\Omega_0)-\mathcal I_T^{\,0}(B;\Omega_0)}{S(T)}=0 .
\]
For a bounded Borel region \(\Omega\), write \(\mathcal R_s\) for multiplication of each exceptional block by \(s_{jl}\in[1/2,2]\), leaving all other blocks fixed. In frames adapted to the flag the enlarged region is
\[
 \Omega^+=\overline{\{\,k\,\mathcal R_s(C)\,k^{-1}:
    C\text{ block upper triangular},\ kCk^{-1}\in \Omega,\ k\in\textstyle\prod_v\mathcal K_v\,\}} .
\]
For any bounded Borel range \(\Omega_1\) of height coordinates one has
\[
 \mathcal I_T(\mathbf1_{\Omega};\Omega_1)\ \le C\,\mathcal I_T^{\,0}(\mathbf1_{\Omega^+};\Omega_1).
\]
The domination holds before integration in the remaining orbital shapes.
\end{lem}
\begin{proof}
\stage{Independent ratios}
Fix the orbit coordinates, and replace the unit groups in each run by a common finite-index subgroup $U$, using probability on the finite covers. Write its independent centralizer variables as $h_j\in E_\infty^1/U$. In an ordinary upper entry the action $C_{jl}\mapsto h_j C_{jl}h_l^{-1}$ has real Jacobian one. In \eqref{exc} it preserves $W$, with Jacobian one there too: on the projected space it acts by multiplication by $h_j/h_l$, of norm one. Hence after changes of the integration variables the only remaining dependence is through the successive ratios on adjacent equal-primary edges. Writing \(r_j=h_j/h_{j+1}\) gives the mutually inverse Haar-preserving maps
\[
 (h_1,\ldots,h_e)\longleftrightarrow(r_1,\ldots,r_{e-1},h_e),
 \qquad h_j=r_j\cdots r_{e-1}h_e,\qquad
 \prod_{j=1}^e dh_j=dh_e\prod_{j<e}dr_j .
\]
Thus the ratios are independent after the ordinary entries and the \(W\)-coordinates have been integrated, with the orbital shapes still fixed.
\stage{An upper bound before taking limits}
Use independent real dilations \(s_e\in[1,1+\varepsilon]\) on the entire physical block at each exceptional edge \(e\), with \(ds_e/(s_e\log(1+\varepsilon))\). Write \(\mathcal I_{T,\varepsilon}\) for the integral with these dilations inserted. They commute with the block shape maps. In the decomposition \(V_e\oplus W_e\), put \(m_e=[E_e:\Q]\), \(w_e=\dim_\R W_e\). Changing variable in \(W_e\) multiplies the measure by \(s_e^{-w_e}\). Relative to ordinary entry measure, the resulting density in the other coordinate is
\[
 r_{e,\varepsilon}(z)=\covol(I_e)\rho_{e,\varepsilon}^{\,W}(z),
 \qquad
 (1+\varepsilon)^{-w_e}\rho_{I_e}(z)
  \le\rho_{e,\varepsilon}^{\,W}(z)\le\rho_{I_e}(z),
\]
using the density of Lemma~\ref{shell} with this \(\varepsilon\). Here \(\covol\Lambda_e=\covol(\Lambda_e\cap W_e)\covol(I_e)\).
For \(\varepsilon=1\) the product over edges is bounded by a model constant. For \(Y\in \Omega\) every dilated image is in \(\Omega^+\). Thus for every height range \(\Omega_1\)
\[
 \mathcal I_T(\mathbf1_{\Omega};\Omega_1)
 \le\mathcal I_{T,1}(\mathbf1_{\Omega^+};\Omega_1)
 \le C\,\mathcal I_T^0(\mathbf1_{\Omega^+};\Omega_1).
\]
These inequalities apply with the other orbital shapes fixed. Applying Lemma~\ref{continuouslemma} to a ball containing \(\Omega^+\) now gives the bound \(O_M(|\Omega_1|S(T))\) for \(\Omega\subseteq\{\|Y\|\le M\}\), uniformly in the measurable height range. The enlargement is essential even in the standard block for \(K=\Q,\ f=t^2\), where the original measure \(\sum_{a\in\Z\setminus\{0\}}\delta_a\) charges the Lebesgue-null set \(\{1\}\).

\stage{Uniformity in the shapes}
For the asymptotic one must guarantee that every projected norm is large in the original field coordinates, even for unbounded orbital shapes. In a fixed integral model the absolute real determinant of the two-layer map $p_i(X)$ is
\[
 c_j |Nz|\,(H_j/H_{j-1})/(H_{j+1}/H_j),\qquad c_j>0 .
\]
This identity is a determinant identity for an arbitrary real upper block: its commutator part is killed and $p_i'(A_j)$ on the intertwiner part is invertible. It continues to hold on all the continuously filled fibers. On a regular compact support the same physical determinant is comparable to $T^{qd^2}$, as in \eqref{index}. For nonempty \(\Omega_0\) and \(\mathcal E\), set
\[
 2\gamma=\min_{\substack{x\in\Omega_0\\(j,j+1)\in\mathcal E}}
                  (x_{j+1}-x_j)>0.
\]
The two determinant formulas imply, uniformly in the orbital shapes,
\[
 \log|Nz|=(x_{j+1}-x_j)\log T+O(1)
       \ge2\gamma\log T-O(1),\qquad |Nz|\gg T^{2\gamma}.
\]
If \(\Omega_0\) is empty both integrals vanish. If \(\mathcal E\) is empty their entry measures coincide. The determinant calculation applies to every ordinary block in the support of the smoothed physical test as well. The transverse Jacobian also identifies the limiting density for fixed \(\varepsilon\). Suppress the edge index. Put \(t=|Nz|\), \(\alpha=(1+\varepsilon)^{-m}\) and \(A_I(y)=\#\{[z']\in(I\setminus0)/U:|Nz'|\le y\}\). An orbit of norm \(y\) contributing at \(t\) has \(s=(t/y)^{1/m}\). In the shell calculation of Lemma~\ref{shell} its extra weight is \(s^{-w}=(y/t)^{w/m}\). Away from the null set of endpoints,
\[
 \frac{\rho^{W}_{e,\varepsilon}(z)}{\rho_I(z)}
  =\frac{\displaystyle\int_{(\alpha t,t]}(y/t)^{w/m}\,dA_I(y)}
          {A_I(t)-A_I(\alpha t)}
 \ \longrightarrow\ \frac{\displaystyle\int_\alpha^1 u^{w/m}\,du}{1-\alpha}.
\]
Indeed \(A_I(y)=c_Iy+o(y)\) implies uniform convergence \(A_I(tu)/t\to c_Iu\) on \([\alpha,1]\).
Writing \(r=w/m\), partial summation yields
\[
 \begin{aligned}
 t^{-1}\!\int_{(\alpha t,t]}(y/t)^r\,dA_I(y)
 &=\frac{A_I(t)}t-\alpha^r\frac{A_I(\alpha t)}t
        -r\int_\alpha^1u^{r-1}\frac{A_I(tu)}t\,du\\
 &\longrightarrow c_I\left(1-\alpha^{r+1}
                   -r\int_\alpha^1u^r\,du\right)
      =c_I\int_\alpha^1u^r\,du .
 \end{aligned}
\]
The integral term is zero when \(r=0\). Multiplying by \(b_{m,\varepsilon}\) from Lemma~\ref{shell} gives
\[
 r_{e,\varepsilon}(z)
   = b_{m+w,\varepsilon}+o(1),\qquad
 b_{k,\varepsilon}=\frac{1-(1+\varepsilon)^{-k}}{k\log(1+\varepsilon)} .
\]
The error is uniform almost everywhere on the support under consideration since \(t\gg T^\gamma\). With \(b_{\mathfrak m,\varepsilon}=\prod_{e\in\mathcal E} b_{m_e+w_e,\varepsilon}\), the product density and \eqref{continuousfactor} therefore give, for fixed \(\varepsilon,\Omega_0\),
\[
 \mathcal I_{T,\varepsilon}(B;\Omega_0)
       =b_{\mathfrak m,\varepsilon}\mathcal I_T^0(B;\Omega_0)+o(S(T)).
\]

\stage{Removing the smoothing}
Let the support of \(B\) lie in \(\|Y\|\le M\), for the product Frobenius norm, and let \(\omega_B(t)=\sup_{\|Y-Y'\|\le t}|B(Y)-B(Y')|\). In an orthonormal flag frame, distinct blocks are orthogonal for this norm. For \(0<\varepsilon\le1\), if \(B(Y)\) or its dilated value is nonzero, then \(\|Y\|\le M\) and the physical displacement is at most \(M\varepsilon\). The rough bound above therefore gives
\[
 S(T)^{-1}|\mathcal I_T(B;\Omega_0)-\mathcal I_{T,\varepsilon}(B;\Omega_0)|
 \le C_{M,\Omega_0}\,\omega_B(M\varepsilon).
\]
Combine this inequality with the fixed-width formula to obtain
\[
 \limsup_{T\to\infty}
 \frac{|\mathcal I_T(B;\Omega_0)-\mathcal I_T^0(B;\Omega_0)|}{S(T)}
 \le C_{M,\Omega_0}
 \big[\omega_B(M\varepsilon)+\|B\|_\infty|1-b_{\mathfrak m,\varepsilon}|\big].
\]
Each \(b_{k,\varepsilon}\to1\) as \(\varepsilon\downarrow0\), so the right side tends to zero. The same estimates allow a bounded measurable factor of \(x\) on \(\Omega_0\): fixed block dilations preserve \(x\), and the error bounds acquire only the supremum of that factor.
\end{proof}

\subsection{The joint measure after averaging the ambient lattice}
We now retain the ambient lattice as a variable. This permits localization at the free lattice in Section~\ref{fixed}.
\begin{defn}[Truncated counting measure]\label{def:measure}
For \(T>1,\ c>0\) define a measure on \(\mathcal X\times M_n(K_\infty)\times\mathcal H\) by
\[
 \begin{aligned}
 P^f_{T,c}(\Phi)
 &=\frac{1}{S(T)}\int_{\mathcal X}
     \sum_{\substack{X\in\End_\Oa(L)\\\min_X=\chi_X=f}}
     \mathbf1_{\{h_F(L,X)\ge c\ \forall\,0<F<A\}}\\[-2pt]
 &\hspace{38mm}\cdot\Phi\big(L,X/T,x_T(L,X)\big)\,d\mu(L).
 \end{aligned}
\]
Here \(x_T(L,X)\) is given by the formula in Section~\ref{introflag} with \(L\) as ambient lattice. Omitting the third argument of a test denotes the \((L,Y)\)-marginal, \(Y=X/T\). The measure for \(g\) uses its separated ordering.
\end{defn}
\begin{prop}[The averaged distribution]\label{avg}
Put \(P_f=P,\ P_g=Z\), and fix \(c>0\) sufficiently small for localization at \(\Oa^n\). There is one constant \(A_0>0\) such that for \(\xi\in\{f,g\}\), bounded continuous \(\psi\), and continuous \(B\) of regular compact support,
\begin{equation}\label{marg}
 \lim_{T\to\infty}P^\xi_{T,c}\big(\psi(x_T)B(Y)\big)
       =A_0\!\left(\int_{P_\xi}\psi(x)\,dx\right)\eta_n(B).
\end{equation}
In particular, the marginal limit is \(A_0\vol(P_\xi)\eta_n\) on regular tests.
\end{prop}
\begin{proof}
\stage{Equality of the two continuous constants}
Fix a tuple of projective layers in Lemma~\ref{continuouslemma}. The bijection of diagonal models is \(A_j\mapsto A_j+a_{ij}1\). Their integral centralizers, entry lattices, and Haar measures coincide. In the coordinates of \eqref{continuousfactor} the two angular matrices differ by
\[
 T^{-1}k\,\operatorname{diag}(a_{ij}1_{b_j})_j\,k^{-1}.
\]
Thus their integrals of \(B\) differ by at most \(C_M\omega_B(C/T)\), by the uniform ball bound there. Both have limit \(c_{\mathfrak m}\eta_n(B)\). Summing over the finite models gives \(A_0=\sum_{\mathfrak m}c_{\mathfrak m}>0\). The companion model on free layers ensures that the sum is nonempty.

\stage{Recovering the whole chamber}
Fix \(0<\delta<1\), and take \(K_{\xi,\delta}=\{x\in P_\xi:
\operatorname{dist}(x,\partial_{\mathcal H} P_\xi)\ge\delta\}\).
On this interior range Corollary~\ref{chamber} removes the proper-ideal cuts on the original configurations. Lemma~\ref{joint} then replaces their edge measures, also in the presence of \(\psi(x)\); in the separated case they are ordinary from the outset. By the exact height factorization \eqref{continuousfactor}, the normalized contribution of each model tends to
\[
 c_{\mathfrak m}\eta_n(B)\int_{K_{\xi,\delta}}\psi(x)\,dx .
\]
For the remaining original configurations with \(B(Y)\ne0\), the same corollary gives \(\operatorname{dist}(x,P_\xi)\le C/\log T\). Remove the cuts for an upper bound. The rough inequality of Lemma~\ref{joint} and \eqref{continuousfactor} bound the mass per unit \(dx\) by a constant on the enlarged ball. Let \(P_\xi^{+h}=\{x\in\mathcal H:\operatorname{dist}(x,P_\xi)\le h\}\). Its finitely many faces have bounded area in a fixed neighborhood, so, for \(\ell>1\),
\[
 \big|P_\xi^{+C/\log T}\setminus K_{\xi,\delta}\big|
     \ll_{P_\xi}\delta+1/\log T .
\] Consequently
\[
 \limsup_{T\to\infty}\left|
 P^\xi_{T,c}(\psi B)-A_0\eta_n(B)\int_{P_\xi}\psi\,dx\right|
 \ll_{B,K,f,g}\|\psi\|_\infty\delta .
\]
Let \(\delta\downarrow0\). If \(\ell=1\), use the entire one-point chamber from the start.
\end{proof}

\section{From a random lattice to the free lattice}\label{fixed}
Proposition~\ref{avg} determines the height and matrix marginal of the measure in Definition~\ref{def:measure}. To recover the value at \(\Oa^n\), we keep the lattice as a coordinate of the measure. A short vector detects escape from compact sets. A vector near the kernel of \(Y_v^{\,n-1}\) detects the proper homogeneous measures that could survive in the limit.

\subsection{A uniform bound for vectors generating an ideal}
Fix an integral conjugacy class of cyclic matrices with polynomial \(f\). In the model \(A=K[t]/(f)\), a representative gives a full \(R=\Oa[t]/(f)\)-lattice \(M_*\). Fix the reference frame identifying it with the free lattice in physical coordinates. For \(0<F\lhd A\), define
\[
 M_F=M_*\cap F,\quad {\rm Gen}_{M_*}(F)=\{w\in M_F:Aw=F\},\quad
 h_F(g)=\covol_{gF_\infty}(gM_F),\quad N=qn .
\]
With \(J=\operatorname{rad} A\), the radical layer is \(W=JF\) and the top quotient is \(F/JF\).
In particular \(h_A(g)=\covol(\Oa^n)\) since \(|\det_\R g|=1\) in the reference frame.

\begin{lem}[Class-uniform density for marked generators]\label{markdensity}
For every class lattice \(M_*\) there is a probability average over multiplication operators \(\xi\in A_\infty^\times\), depending on the class but independent of \(F\) and \(g\), with the following property. For every ideal \(0<F\le A\), every frame \(g\in G\), and every Borel \(\Omega\subset gF_\infty\),
\begin{equation}\label{generatorbound}
 \mathbb E_\xi\!\sum_{w\in{\rm Gen}_{M_*}(F)}
       {\bf1}_{\Omega}(g\xi w)
       \ \le\ \frac{C(K,f)}{h_F(g)}\vol_{gF_\infty}(\Omega).
\end{equation}
If \(d_\xi=|\det_\R\xi|^{1/N}\) and \(\widehat\xi=\xi/d_\xi\), one may arrange simultaneously
\[
   1\le d_\xi\le2,\qquad
   2^{-N}\le \frac{h_F(g\widehat\xi)}{h_F(g)}\le 2^N
           \quad(0<F\le A).
\]
\end{lem}

\begin{proof}
\stage{A quotient with one field coordinate per active primary}
Write \(A=\prod_i E_i[\epsilon_i]/(\epsilon_i^{e_i})\) using a separable \(K\)-algebra section of the semisimple quotient, and let \(J\) be the radical. For \(F\) set
\[
 W=JF,\quad V=F/W,\quad M_W=M_*\cap W,\quad
 \overline M=M_F/M_W.
\]
With the induced subspace and quotient volumes,
\begin{equation}\label{covolidentity}
 h_F(g)=\covol(gM_W)\covol(\overline{gM_F}).
\end{equation}
Indeed \(F=\prod_i\epsilon_i^{k_i}A_i\) for unique \(0\le k_i\le e_i\), where \(A_i=E_i[\epsilon_i]/(\epsilon_i^{e_i})\). Put \(\ell_i=e_i-k_i\) and \(\mathcal J_F=\{i:\ell_i>0\}\). Taking the coefficient of \(\epsilon_i^{k_i}\) gives
\[
 F/JF\simeq\prod_{i\in \mathcal J_F} E_i,\qquad
 Aw=F\ \Longleftrightarrow\
 [\epsilon_i^{k_i}]w_i\ne0\quad(i\in \mathcal J_F).
\]
Indeed \(w_i=\epsilon_i^{k_i}u_i\) generates the \(i\)th ideal exactly when the constant coefficient of \(u_i\) is a unit of the field \(E_i\).

\stage{Exact conditional density along \(W\)}
Average over the unipotent group \(1+J_\infty\) modulo its periods preserving \(M_*\), with probability Haar measure. It is a compact quotient: the logarithm identifies this abelian group with \(J_\infty\), and exponentials of a sufficiently divisible full integral lattice in the radical of \(R\) preserve \(M_*\), together with their inverses. For a fixed generating top vector the lifts form a coset of \(M_W\). If \(w\) is one lift, \(J_\infty w=W_\infty\). Therefore \(1+J_\infty\) acts transitively on the affine fiber, preserving its Lebesgue measure. The averaged counting measure there is a multiple of Lebesgue measure.

To compute the multiple, let \(k=\dim_\R W_\infty>0\), \(h_W=\covol(gM_W)\), and choose bounded representatives of the compact unipotent quotient for this class. Their affine lattices have fundamental cells of some common diameter \(D_*\), and all have covolume \(h_W\) since \(u|_{W_\infty}\) is unipotent. If \(n_u(R)\) counts a ball of radius \(R>D_*\) in the affine fiber, the fundamental-cell inclusions give
\[
 v_k(R-D_*)^k\le h_W n_u(R)\le v_k(R+D_*)^k,\qquad
 \lim_{R\to\infty}\frac{\mathbb E_u n_u(R)}{v_kR^k}=h_W^{-1},
\]
where \(v_k\) is unit ball volume. Transitivity makes this the density on the entire fiber. If \(k=0\), the fiber is a singleton of mass one. In particular, for a generating top vector \(\bar w_0\) and a lift \(w_0\),
\[
 \int_{\mathcal U_*}\sum_{w\in w_0+M_W}\phi(guw)\,du
 =\frac{1}{\covol(gM_W)}
     \int_{g w_0+gW_\infty}\phi(y)\,dy ,
\]
where \(\mathcal U_*=(1+J_\infty)/\operatorname{Stab}_{1+J_\infty}(M_*)\) has mass one. This formula is applied before averaging over the field coordinates.

\stage{Uniform upper density on \(V\)}
In the field coordinates write \(I_i\) for the projection of \(\overline M\) to the \(i\)th active factor. There exists an integer \(b\ge1\), depending on \(R\) alone, with \(b\prod_i\Oa_{E_i}\subset R\) in the chosen semisimple section: clear the denominators of an integral basis once. Applying the multiples of the primary idempotents to \(M_F\) gives
\[
 b\prod_{i\ {\rm active}} I_i\ \subseteq\ \overline M\
     \subseteq\ \prod_{i\ {\rm active}} I_i .
\]
Each \(I_i\) is stable under the fixed order generated by the image of \(t\) on that layer. After scaling, the maximal-order saturation has one of finitely many ideal-class models. The conductor then leaves only finitely many possibilities for \(I_i\). For \(u_i\in\Oa_{E_i}^\times\), \(u_i\equiv1\bmod b\Oa_{E_i}\), the embedded multiplier and its inverse belong to \(R\):
\[
 (1,\ldots,u_i,\ldots,1)=1+\varepsilon^{(i)}(u_i-1)\in R,\qquad
 (1,\ldots,u_i^{-1},\ldots,1)\in R .
\]
Here \(\varepsilon^{(i)}\) is the primary idempotent. Thus each multiplier preserves every \(R\)-lattice \(M_*\). Let \(U_i\) be these congruence-unit groups. Use bounded measurable representatives of each compact quotient and take the single product probability
\[
 \xi=u\,(s_i h_i)_i,\qquad
 du\prod_i \left(dh_i\,\frac{ds_i}{s_i\log2}\right),\qquad
 u\in\mathcal U_*,\quad h_i\in E_{i,\infty}^{1}/U_i,\quad 1\le s_i\le2 .
\]
All factors commute. A change of representative permutes the vectors being summed. This law uses the full primary list and is chosen before \(g\) or \(F\).

To bound the sum over generating top vectors, enlarge it to the product of the sums over \(I_i\setminus0\). Lemma~\ref{shell} bounds the product density by a constant times \(1/\covol(\prod I_i)\) in model coordinates. Since the index of \(\overline M\) in this product is at most \(b^{\dim_\R V_\infty}\), the bound is at most \(C/\covol(\overline M)\). Under any change to a physical quotient frame, reciprocal covolume and density acquire the same Jacobian. On \(W_\infty\) the norm-one multipliers have determinant of absolute value one; the \(s_i\) have determinant at least one. Thus the conditional density along \(W\) remains at most \(1/\covol(gM_W)\). Let \(\pi:gF_\infty\to gF_\infty/gW_\infty\) be the quotient map and use quotient Lebesgue measure \(d\bar y\). For every nonnegative Borel test \(\phi\) the two density bounds give
\[
 \mathbb E_\xi\sum_{Aw=F}\phi(g\xi w)
 \le \frac{C}{\covol(\overline{gM_F})\,\covol(gM_W)}
       \int\!\int_{\pi^{-1}(\bar y)}\phi(y)\,dy_W\,d\bar y
 =\frac{C}{h_F(g)}\int_{gF_\infty}\phi(y)\,dy .
\]
Here the sum runs over \(w\in M_F\). Taking \(\phi=\mathbf1_{\Omega}\) gives \eqref{generatorbound}. For clarity, let \(C_i\ge1\) bound \(\covol(I_i)\|\rho_{I_i}\|_\infty\) over the finite list of homothety classes at primary \(i\) for the fixed smoothing interval \([1,2]\). The preceding proof permits
\[
 C(K,f)=b^{N}\prod_i C_i,\qquad N=q n.
\]
Indeed, the top index is \(\le b^{\dim_\R V_\infty}\le b^N\), and the exact conditional density on \(W\) contributes no class-dependent factor. The compact set of unipotent representatives used to establish that identity may vary with \(M_*\). Only the resulting identity is used here.

Finally put \(m_i=[E_i:\Q]\), \(N=\sum_i m_i e_i\) and \(k_F=\sum_i m_i \ell_i\). The unipotent and norm-one factors have absolute determinant one on every ideal. Hence
\[
 d_\xi=\prod_i s_i^{m_i e_i/N}\in[1,2],\qquad
 \frac{h_F(g\widehat\xi)}{h_F(g)}
   =|\det_\R(\widehat\xi|F_\infty)|
   =\prod_i(s_i/d_\xi)^{m_i \ell_i}\in[2^{-k_F},2^{k_F}].
\]
This proves the asserted simultaneous bounds since \(k_F\le N\).
\end{proof}

\subsection{Short vectors and vectors near a kernel}
\begin{lem}[Nondivergence and kernel avoidance]\label{mark}
Fix a regular compact matrix set \(\mathcal C\), \(\gamma>0\), and \(R_0>0\). In \(P^f_{T,c}\) restrict to \(Y=X/T\in\mathcal C\) and \(h_F\ge T^\gamma\) for every proper nonzero ideal. Denote integration on this restricted measure by \(\int_*\). Uniformly for \(0<R,\tau\le1\) and \(1\le r<n\),
\begin{align*}
 \int_* \#\{0\ne w\in L:\|w\|<R\}&\ll R^N+T^{-\gamma},\\
 \int_* \#\{0\ne w\in L:\|w\|<R_0,\ 
            \|Y_v^r w_v\|<\tau\ \text{for every }v\}
       &\ll_{R_0}\tau^{q(n-r)}+T^{-\gamma}.
\end{align*}
Constants can depend on \(K,f,\mathcal C,c,\gamma\).
\end{lem}
\begin{proof}
Unfold by integral conjugacy classes and fix \(X_*\) on \(M_*\). Denote its quotient measure on \(G/\Gamma_{X_*}\) by \(d\omega\). Set \(Y(g)=gX_*g^{-1}/T\) and let \(\mathcal A_H\) impose \(Y(g)\in\mathcal C\) and \(h_F(g)\ge H\) for all proper nonzero ideals. The normalized multiplier \(\widehat\xi\) of Lemma~\ref{markdensity} centralizes \(\Gamma_{X_*}\), so right translation by \(\widehat\xi\) preserves \(d\omega\). For either vector region \(\mathcal V(Y)\) in the statement,
\[
 \begin{aligned}
 Y(g\widehat\xi)&=Y(g),&
 g\widehat\xi\in\mathcal A_H&\ \Longrightarrow\ g\in\mathcal A_{2^{-N}H},\\
 g\widehat\xi w\in \mathcal V(Y)&\ \Longrightarrow\
 g\xi w=d_\xi g\widehat\xi w\in 2\mathcal V(Y) &&(1\le d_\xi\le2).
 \end{aligned}
\]
Change variables and average these inequalities to obtain
\[
 \int_{\mathcal A_H}\sum_{0\ne w\in M_*}\mathbf1_{\mathcal V(Y)}(gw)\,d\omega
 \le \int_{\mathcal A_{2^{-N}H}}\mathbb E_\xi
             \sum_{0\ne w\in M_*}\mathbf1_{2\mathcal V(Y)}(g\xi w)\,d\omega .
\]
Partition the right-hand sum by \(F=Aw\) and take \(H=T^\gamma\). For \(F=A\) its volume factor is \(O(R^N)\) in the ball case. In the tube case write the local singular values of \(Y_v^r\) in decreasing order. Every limit of contributing matrices in \(\mathcal C\) as \(T\to\infty\) is regular nilpotent, so \(\sigma_{n-r}(Y_v^r)\ge c_{\mathcal C}>0\) on the contributing support for all large \(T\). In right singular coordinates \(u\) this gives
\[
 \|Y_v^r u\|^2=\sum_j \sigma_j^2 |u_j|^2<4\tau^2
 \ \Longrightarrow\ |u_j|\le 2\tau/c_{\mathcal C}\quad(1\le j\le n-r).
\]
The other coordinates are bounded by \(2R_0\), giving total volume
\(O_{R_0,\mathcal C}(\prod_v\tau^{q_v(n-r)})=O_{R_0,\mathcal C}(\tau^{q(n-r)})\).

To see the uniformity in the sum over classes, let \(W_{\mathfrak a}(H')\) denote the unmarked quotient mass on \(\mathcal A_{H'}\) for the class \(\mathfrak a\). For a fixed \(g\) in the relaxed region put
\(\mathcal V_F=\vol_{gF_\infty}(2\mathcal V(Y)\cap gF_\infty)\).
The ideals are
\[
 F=\prod_i\epsilon_i^{k_i}E_i[\epsilon_i]/(\epsilon_i^{e_i}),
 \qquad 0\le k_i\le e_i,\qquad
 \#\{F\lhd A\}=\prod_i(e_i+1).
\] Lemma~\ref{markdensity} gives the same bound on every class,
\[
 C(K,f)\sum_{0<F\le A}\frac{\mathcal V_F}{h_F(g)}
 \ \ll\ \mathcal V_A+T^{-\gamma}.
\]
Here \(h_A\) is fixed, and all proper slices have bounded volume since they lie in a fixed ball in dimension at most \(N\). Put \(v_{\mathcal V}=R^N\) for the ball or \(v_{\mathcal V}=\tau^{q(n-r)}\) for the tube, so that \(\mathcal V_A\ll v_{\mathcal V}\). For a continuous regular outer test \(B_0\ge\mathbf1_{\mathcal C}\) we obtain, by nonnegativity and unfolding,
\[
 \begin{aligned}
 \int_*\#\{0\ne w\in L:w\in \mathcal V(Y)\}
 &\ll (v_{\mathcal V}+T^{-\gamma})\,S(T)^{-1}
                 \sum_{\mathfrak a} W_{\mathfrak a}(2^{-N}T^\gamma)\\
 &\le C(v_{\mathcal V}+T^{-\gamma})P^f_{T,c}(B_0)
 \ll v_{\mathcal V}+T^{-\gamma}.
 \end{aligned}
\]
In the second line \(2^{-N}T^\gamma\ge c\) for large \(T\), and Proposition~\ref{avg} bounds the last unmarked mass. This sums the classes before taking a limit in \(T\).
\end{proof}

\begin{lem}[Tightness and kernel avoidance on regular support]\label{regularlimits}
For every nonnegative \(B\in C_c\) of regular compact support, the measures \(B(Y)P^f_{T,c}\) have uniformly bounded mass and are tight as \(T\to\infty\). Each subsequential weak limit is supported over \(\mathcal N^{\rm reg}\times P\) and satisfies
\begin{equation}\label{avoid}
 (L\setminus\{0\})\cap\prod_v\ker Y_v^{n-1}=\varnothing
 \quad\text{almost everywhere}.
\end{equation}
\end{lem}
\begin{proof}
Proposition~\ref{avg} bounds the masses. Corollary~\ref{chamber} confines the height coordinates to a fixed compact neighborhood of \(P\); every matrix limit is nilpotent and lies on the regular support.
Discarding \(h_F<T^\gamma\) for a proper ideal costs \(O_B(\gamma+1/\log T)\) by \eqref{cuts} and the wall-strip bound. Existence of a lattice vector in a strict ball and a strict tube defines an open condition. Markov's inequality, the strip estimate and Lemma~\ref{mark} give, for each \(\gamma>0\),
\[
 \begin{aligned}
 \limsup_{T\to\infty}P^f_{T,c}(\text{some }0<\|w\|<R)
     &\ll C_1\gamma+C_\gamma R^N,\\
 \limsup_{T\to\infty}P^f_{T,c}(\text{some }0<\|w\|<R_0,\
                        \max_v\|Y_v^{\,n-1}w_v\|<\tau)
     &\ll C_1\gamma+C_{\gamma,R_0}\tau^q .
 \end{aligned}
\]
Both estimates are for the weighted measure \(B(Y)P^f_{T,c}\). Portmanteau passes them to any weak limit on regular support. For each fixed \(R_0\), first let \(\tau\downarrow0\) at fixed \(\gamma\), then let \(\gamma\downarrow0\). The tube probability in the limit is zero for that \(R_0\). Taking the countable union over \(R_0=1,2,\ldots\) gives \eqref{avoid}. For tightness the first display gives
\[
 \lim_{R\downarrow0}\limsup_{T\to\infty}
 P^f_{T,c}(\exists\,w\in L:\ 0<\|w\|<R)
 \le C_1\gamma\qquad(\gamma>0)
\]
on the fixed regular test. Now let \(\gamma\downarrow0\) and apply Mahler's criterion.

The real Mahler criterion applies here as follows. In a fixed integral basis of $\Oa_K^n$ let $g_j B_j$ be convergent real frames with $g_j\in G$ and $B_j\in\GL_N(\Z)$. For each generator $a$ of $\Oa_K$ the matrices $B_j^{-1}aB_j=(g_j B_j)^{-1}a(g_j B_j)$ are integral and bounded. On a subsequence they are constant for all generators. Then $B_j B_{j_0}^{-1}$ is integral, invertible and $\Oa_K$-linear. This gives convergent representatives modulo $\Gamma$ and proves the required properness.

\end{proof}

\subsection{Unipotent invariance and the determinant}
Put \(K_\infty^1=\{z\in K_\infty^\times:\prod_v|z_v|^{q_v}=1\}\),
\(G'=\prod_v\SL_n(K_v)\), and define
\[
 \mathcal T_{\det}=K_\infty^1/\Oa^\times,\qquad
 \delta(g\Gamma)=[(\det g_v)_v],\qquad \mathcal X_d=\delta^{-1}(d).
\]
The compact group \(\mathcal T_{\det}\) has probability Haar measure. Write \(m_d\) for the \(G'\)-invariant probability on the translated fiber \(\mathcal X_d\). For \(Y\in\mathcal N^{\rm reg}\) set
\[
 U_Y=\prod_v\left\{\exp\!\left(\sum_{j=1}^{n-1}t_{v,j}Y_v^j\right):
                               t_{v,j}\in K_v\right\},\qquad
 K(Y)=\prod_v\ker(Y_v^{n-1}).
\]
The unit theorem gives compactness of \(\mathcal T_{\det}\). The fiber over the identity is \(G'/\SL_n(\Oa)\), since \(\det\Gamma=\Oa^\times\); translation identifies the other fibers with it. 
\begin{lem}[Unipotent invariance of the conditional measure]\label{flowinvariance}
Let \(\nu\) be a subsequential limit as in Lemma~\ref{regularlimits}. Disintegrate it over \((Y,x,d)\), \(d=\delta(L)\). For almost every base point its conditional lattice probability satisfies
\[
        u_*\nu_{Y,x,d}=\nu_{Y,x,d}\qquad(u\in U_Y).
\]
\end{lem}
\begin{proof}
Fix an infinite place \(v\), \(1\le j<n\), \(t\in K_v\), and a class representative \(X_*\). Put, at \(v\),
\[
 Z_{j,T}=(X_{*,v}/T)^j-
          \frac{\operatorname{tr}((X_{*,v}/T)^j)}{n}I_n,\qquad
 u_{j,T}(t)=\exp(tZ_{j,T}),
\]
and put \(u_{j,T}(t)=1\) elsewhere. This translates the unfolded class on the right, preserves Haar measure and has determinant one. Its action on a physical lattice is on the left by \(gu_{j,T}(t)g^{-1}\). Each invariant ideal is preserved, and
\[
 \log\frac{h_F(gu_{j,T}(t))}{h_F(g)}
 =q_v\operatorname{Re}\big(t\,\operatorname{tr}(Z_{j,T}|F_v)\big)
 =O_{t,K,f}(T^{-j}).
\]
The traces involve only the fixed eigenvalues with multiplicities. In physical coordinates the local flow is \(\exp(t[(gX_*g^{-1}/T)_v^j-\operatorname{tr}((X_{*,v}/T)^j)I_n/n])\). The trace depends only on \(f\). Thus convergence to \(\exp(tY_v^j)\) along \(gX_*g^{-1}/T\to Y\) is uniform on fixed compact matrix sets, simultaneously in the integral class.

The prefix trace calculation also gives
\[
 \|x_T(gu_{j,T}(t))-x_T(g)\|_\infty
        \ll_{t,K,f}T^{-j}/\log T .
\]
A cutoff crossing has \(h_F/c=\exp(O(T^{-j}))\). By \eqref{cuts} it lies within \(O(1/\log T)\) of \(x(I)+Q(b(I))=0\), whose mass is \(o(1)\) by the rough strip bound. The translations fix \(Y\), the determinant coset \(d\), and hence \(B(Y)\). On the compact height range from Corollary~\ref{chamber}, the height displacement tends to zero uniformly.

For \(\Phi\in C_c\) of \((L,Y,x,d)\), change variables in the unfolded integral and omit the cutoff crossings. Uniform convergence of the physical flows and the height displacement gives
\[
 \int(\Phi\circ\Theta_{v,j,t}-\Phi)\,d\nu=0,\qquad
 \Theta_{v,j,t}(L,Y,x,d)=(\exp(tY_v^j)L,Y,x,d).
\]
Since \(\Theta_{v,j,t}\) fixes the base \((Y,x,d)\), disintegration gives
\[
 (\exp(tY_v^j))_*\nu_{Y,x,d}=\nu_{Y,x,d}
\]
for almost every base point, simultaneously on a countable dense set of \(t\).
Continuity extends the equality to all \(t\in K_v\).
The local powers generate \(U_Y\), proving the assertion.
\end{proof}

\begin{lem}[Haar measure on a determinant fiber]\label{rigid}
Let \(Y\in\mathcal N^{\rm reg}\), \(d\in\mathcal T_{\det}\), and let \(\nu\) be a probability on \(\mathcal X_d\). Then
\[
 \left.
 \begin{aligned}
    u_*\nu&=\nu &&(u\in U_Y),\\
    \nu\{L:(L\setminus\{0\})\cap K(Y)\ne\varnothing\}&=0
 \end{aligned}\right\}
 \quad\Longrightarrow\quad \nu=m_d .
\]
\end{lem}
\begin{proof}
\stage{The algebraic group attached to an ergodic component}
Ratner's measure theorem \cite{Ratner} gives an ergodic component as invariant probability \(m\) on \(HL_*\simeq H/\Delta\), with \(U_Y\subset H\subseteq G'\) and \(\Delta=\operatorname{Stab}_H(L_*)\) a lattice. We may use connected \(H\). Indeed \(U_Y\subset H^\circ\), and the \(H^\circ\)-orbits in \(H/\Delta\) are disjoint open invariant sets. Full support and ergodicity leave only one. Replace \(H\) by \(H^\circ\) and \(\Delta\) by its intersection with \(H^\circ\).

Use the rational frame of \(L_*\) and set \(\mathbf S=\overline\Delta^{\,\rm Zar}\) in the restriction of scalars of \(\SL_n\). The elements of \(\Delta\) are rational in this frame, so \(\mathbf S\) is defined over \(\Q\). By Chevalley's stabilizer-of-a-line theorem \cite[Theorem 4.27]{Milne}, there is an algebraic representation and a projective point \([v]\) with stabilizer \(\mathbf S\). Consider the equivariant map
\[
 \Phi:H/\Delta\longrightarrow\mathbb P(V),\qquad h\Delta\longmapsto h[v].
\]
Its image law \(\lambda=\Phi_*m\) is \(U_Y\)-invariant and ergodic. Take a real one-parameter generator with nilpotent derivative \(\mathsf N\). For a nonzero \(w\) let \(j=\max\{i:\mathsf N^iw\ne0\}\). Then
\[
 [\exp(t\mathsf N)w]=\left[\sum_{i=0}^{j}\frac{t^i}{i!}\mathsf N^iw\right]
       \longrightarrow[\mathsf N^jw]\qquad(t\to+\infty).
\]
For \(\lambda\)-almost every \([w]\), Poincar\'e recurrence for the time-one map supplies integers \(m_k\to\infty\) with
\[
 [w]=\lim_k[\exp(m_k\mathsf N)w]=[\mathsf N^jw].
\]
If \(j>0\), then \(\mathsf N^jw=cw\) for some \(c\ne0\) implies \(0=\mathsf N^{j+1}w=c\mathsf N w\), contradicting \(\mathsf N^jw\ne0\). Thus \(j=0\) and \(\mathsf N w=0\). Intersect the full-measure sets for the finitely many generators. The action of \(U_Y\) on that intersection is trivial. Ergodicity of \(\lambda\) forces \(\lambda=\delta_{[w_0]}\). Full support of \(m\) and continuity of \(\Phi\) then give
\[
 \Phi(h\Delta)=\Phi(\Delta)=[v]\quad(h\in H),\qquad
 \Delta\subset H\subset\mathbf S(\R),\qquad
 \overline H^{\,\rm Zar}=\overline\Delta^{\,\rm Zar}.
\]
Recurrence gives fixed projective points. Ergodicity is the additional step that makes the entire image a single point.
\stage{A reducible projection gives a forbidden vector}
In the rational structure given by \(L_*\), let \(\mathbf H_\sigma\) be the projection of \(\mathbf S_\C\) to the factor for an embedding \(\sigma:K\hookrightarrow\C\). Suppose it preserves an \(r\)-plane, \(0<r<n\). Its fixed-plane locus satisfies
\[
 \varnothing\ne\operatorname{Gr}(r,n)^{\mathbf H_\sigma}
    \subseteq\operatorname{Gr}(r,n)^{\exp(\C Y_\sigma)}
    =\{\ker Y_\sigma^r\}
\]
in the corresponding frame. Indeed the Jordan module is \(\C[t]/(t^n)\), whose unique submodule of dimension \(r\) is \((t^{\,n-r})=\ker t^r\). The projected group is defined over \(\sigma(K)\): the stabilizer of \(\sigma\) in \(\operatorname{Aut}(\C/\Q)\) preserves this projection of the rational group \(\mathbf S\). In Pl\"ucker coordinates \([p_I]\) for the unique invariant plane, choose \(p_{I_0}\ne0\). Uniqueness gives
\[
 \tau(p_I/p_{I_0})=p_I/p_{I_0}
 \quad(\tau\in\operatorname{Aut}(\C/\sigma(K))),
 \qquad p_I/p_{I_0}\in\sigma(K).
\]
On the affine chart \(p_{I_0}\ne0\), these ratios include the signed entries of the row-reduced matrix of the plane. This constructs the descended \(K\)-plane explicitly. At every other embedding the conjugate plane is fixed by the corresponding projected group and consequently equals the local Jordan plane in physical coordinates. If \(V\) denotes the descended plane, choose \(0\ne w\in V\cap L_*\) by clearing denominators. Since \(H\) preserves all these planes,
\[
       hw\in (hL_*\setminus\{0\})\cap\prod_v\ker Y_v^r
                \subseteq(hL_*\setminus\{0\})\cap K(Y)\qquad(h\in H).
\]
Thus every point of this orbit has a forbidden vector, contradicting the almost-everywhere avoidance on the component.
\stage{Irreducible projections generate the whole group}
The Zariski closure of the connected group \(H\) is geometrically connected: \(H\) lies in the identity component of its complex closure and is Zariski dense there. Thus each projected algebraic group is connected. In the remaining case its Lie algebra \(\mathfrak h\) acts irreducibly on \(\C^n\), since the stabilizer of a subspace is closed and contains the connected group whenever its Lie algebra contains \(\mathfrak h\). Let \(R_u\) be the unipotent radical of this projected group. Its fixed subspace \(V^{R_u}\) is nonzero by triangularization. Normality and irreducibility give
\[
 0\ne V^{R_u}\subseteq\C^n\text{ invariant}\ \Longrightarrow\
 V^{R_u}=\C^n\ \Longrightarrow\ R_u=\{1\}.
\]
Schur's lemma places the connected center in the scalar subgroup of \(\SL_n\), which is finite. Its Lie algebra is therefore zero and \(\mathfrak h\) is semisimple. Choose an \(\mathfrak{sl}_2\)-triple \((Y,H,F)\) in \(\mathfrak h\) by Jacobson--Morozov. The regular Jordan block makes \(\C^n\) the irreducible module of dimension \(n\). Clebsch--Gordan gives
\[
 \mathfrak{sl}_n(\C)\simeq\bigoplus_{j=1}^{n-1}\operatorname{Sym}^{2j}(\C^2),
 \qquad [Y,Y^j]=0,\qquad[H,Y^j]=2jY^j.
\]
The one-parameter groups defining \(U_Y\) give \(Y^j\in\mathfrak h\) for \(1\le j<n\). Each \(Y^j\ne0\) is the highest-weight vector of the unique summand of weight \(2j\). In particular
\[
 \mathfrak h\supseteq
 \bigoplus_{j=1}^{n-1}
 \operatorname{span}_{\C}\{(\operatorname{ad}F)^bY^j:0\le b\le2j\},
 \qquad
 \dim\mathfrak h\ge\sum_{j=1}^{n-1}(2j+1)=n^2-1.
\]
Thus every projection is \(\mathfrak{sl}_n(\C)\).

It remains to exclude relations between projections. Write
\(\mathfrak s=\operatorname{Lie}\mathbf S_\C\subseteq
 \bigoplus_\sigma\mathfrak g_\sigma\), \(\mathfrak g_\sigma=\mathfrak{sl}_n(\C)\), and \(\mathfrak i_\sigma=\mathfrak s\cap\mathfrak g_\sigma\). The independent local flows give a nonzero element supported on the factor of each real embedding. At a complex place the real flow derivatives complexify to \(a=(Y,\bar Y)\) and \(b=(iY,-i\bar Y)\); \((a-ib)/2=(Y,0)\) and \((a+ib)/2=(0,\bar Y)\). Thus \(\mathfrak i_\sigma\ne0\) for every embedding. For \(z\in\mathfrak g_\sigma\) choose \(\tilde z\in\mathfrak s\) projecting to \(z\). Then
\[
 [z,w]=[\tilde z,w]\in\mathfrak i_\sigma\quad(w\in\mathfrak i_\sigma)
 \quad\Longrightarrow\quad
 \mathfrak i_\sigma=\mathfrak g_\sigma
 \quad\Longrightarrow\quad \mathfrak s=\bigoplus_\sigma\mathfrak g_\sigma ,
\]
by simplicity. Consequently \(H\) is real Zariski dense in \(G'\). The stabilizer of the real subspace \(\operatorname{Lie} H\) in the adjoint representation is Zariski closed and contains \(H\). It therefore contains \(G'\), so \(\operatorname{Lie} H\) is an ideal. The independent local flow derivatives give a nonzero intersection with each simple real factor. Hence \(\operatorname{Lie}H=\operatorname{Lie}G'\) and, by connectedness, \(H=G'\). Every ergodic component satisfying avoidance is therefore \(m_d\); the hypothesis gives avoidance on almost every component.
\end{proof}

\subsection{The determinant marginal}
\begin{lem}[The ambient lattice on regular support]\label{latticemarginal}
Let \(B\in C_c\) have regular compact support and satisfy the reflection condition of Theorem~\ref{law}. For bounded continuous \(\psi\) and \(\rho\in C_c(\mathcal X)\),
\[
 \lim_{T\to\infty}P^f_{T,c}\big(\rho(L)\psi(x_T)B(Y)\big)
     =\mu(\rho) A_0\eta_n(B)\int_P\psi(x)\,dx .
\]
\end{lem}
\begin{proof}
By linearity take \(B,\psi\ge0\). Lemmas~\ref{regularlimits}, \ref{flowinvariance}, and \ref{rigid} give \(\nu_{Y,x,d}=m_d\) for any subsequential limit of \(B P^f_{T,c}\). It remains to determine the \(d\)-marginal of \(\psi(x)d\nu\).

For \(z\in K_\infty^1\), multiplication of the lattice by \(zI_n\) preserves the matrix and all ideal covolumes:
\[
 \delta(zL)=[z^n]\delta(L),\qquad
 \frac{h_F(zL,X)}{h_F(L,X)}
       =\prod_v|z_v|^{q_v\dim_K F}=1 .
\]
It preserves \(P^f_{T,c}\) by Haar invariance. Thus the weighted determinant marginal \(\omega_\psi\) of \(\nu\) is invariant under \(\mathcal T_{\det}^n\). Local \(n\)th roots show
\[
 \mathcal T_{\det}/\mathcal T_{\det}^n\simeq
 \begin{cases}
 \{1\},&n\text{ odd},\\
 \{\pm1\}^{r_1}/\operatorname{sgn}(\Oa^\times),&n\text{ even}.
 \end{cases}
\]
For even \(n\), a class in the kernel of the sign quotient can be multiplied by a unit to become positive at all real places, after which it has local \(n\)th roots. Their total norm modulus is one. For odd \(n\) the roots exist at every place.

Whenever the sign quotient is nontrivial, each required reflection \(R_v\) preserves \(B\), the prefix heights and the measure, while multiplying \(d\) by the corresponding sign. Therefore \(\omega_\psi\) is invariant under all of \(\mathcal T_{\det}\), and
\[
 \omega_\psi=\nu(\psi)m_{\mathcal T_{\det}},\qquad
 \nu(\rho\psi)=\int_{\mathcal T_{\det}}m_d(\rho)\,d\omega_\psi(d)
             =\nu(\psi)\int_{\mathcal T_{\det}}m_d(\rho)\,dm_{\mathcal T_{\det}}(d)
             =\nu(\psi)\mu(\rho).
\]
Proposition~\ref{avg} gives \(\nu(\psi)=A_0\eta_n(B)\int_P\psi\,dx\). Every subsequence has the same limit for this test, which proves the assertion.
\end{proof}

\begin{example}[An angular test that fails at the free lattice]
Take \(K=\Q(\sqrt3)\), \(s=\sqrt3\), \(f(t)=t^2-st+1\), and let \(v_+,v_-\) be its two real places. We claim
\[
 (X_{v_+})_{21}(X_{v_-})_{21}>0
 \qquad\bigl(X\in M_2(\Oa_K),\ \chi_X=f\bigr).
\]
Indeed, if \(\alpha\) is a root and \(F=K(\alpha)\), then \(\alpha=(s+i)/2\) can be taken to be a primitive twelfth root of unity. The basis \(1,\alpha\) has relative discriminant \(s^2-4=-1\), a unit, so \(\Oa_K[\alpha]=\Oa_F\). The signature of \(F\) is \((0,2)\), and \(|\Delta_F|=|\Delta_K|^2=144\). Minkowski's ideal-class bound is
\[
 \left(\frac4\pi\right)^2\frac{4!}{4^4}\sqrt{144}
 =\frac{18}{\pi^2}<2.
\]
Every ideal class therefore has an integral representative of norm one. Under \(K^2\simeq F\) an \(X\)-stable free lattice becomes a fractional \(\Oa_F\)-ideal. It is principal, and the basis obtained by multiplying \(1,\alpha\) by its generator shows that every such \(X\) is \(hCh^{-1}\), where
\[
 C=\begin{pmatrix}0&-1\\1&s\end{pmatrix},\qquad h\in\GL_2(\Oa_K).
\]
At either real place \(v\) write \(h_v^{-1}e_1=(u,w)^t\). Then
\[
 (X_v)_{21}=\det(e_1,X_ve_1)
   =\det(h_v)\bigl(u^2+s_vuw+w^2\bigr)
   =\det(h_v)\bigl((u+s_vw/2)^2+w^2/4\bigr).
\]
The quadratic factor is positive. Every unit of \(\Z[\sqrt3]\) has norm \(+1\), since \(a^2-3b^2=-1\) is impossible modulo \(3\). Thus the determinants of \(h\) at the two places have the same sign, proving the claim.

The open set \(\{Y:(Y_{v_+})_{21}(Y_{v_-})_{21}<0\}\) contains the nilpotent pair \((E_{21},-E_{21})\), where \(E_{21}\) is the matrix unit. Choose \(0\le B\in C_c\) supported in this set with \(B(E_{21},-E_{21})>0\). The map \((k,U)\mapsto kUk^{-1}\) takes some point to each indicated matrix, so the inverse image of a small neighborhood on which \(B>0\) contains a nonempty open set of positive product measure. Hence
\[
 \eta_2(B)>0,\qquad
 \sum_{\chi_X=f,\ X\in M_2(\Oa_K)} B(X/T)=0\qquad(T>0).
\]
For this \(K,f\), Theorem~\ref{law} fails for the displayed angular test. The entry box satisfies the required reflection symmetry.
\end{example}

\subsection{Localization and removal of the singular set}
\begin{prop}[The cyclic count at a single lattice]\label{localized}
Let \(\psi,B\) satisfy the test conditions of Theorem~\ref{law}, and let \(A_0\) be the constant of Proposition~\ref{avg}. Then
\[
 \lim_{T\to\infty}\frac1{S(T)}
   \sum_{\substack{X\in M_n(\Oa)\\\min_X=\chi_X=f}}
      \psi(x_T(X))B(X/T)=A_0\,\eta_n(B)\int_P\psi(x)\,dx .
\]
The formula also holds for \(B=\mathbf1_{\Omega}\) when \(\Omega\) is bounded Borel, \(\eta_n(\partial \Omega)=0\), and the applicable reflection condition holds.
\end{prop}
\begin{proof}
For \(n=1\) the fiber contains the single matrix prescribed by \(f\), \(x_T=0\), and \(Y\to0\). Assume \(n\ge2\). Write \(M_T(\psi,B)\) for the normalized sum in the statement and \(M_T(B)=M_T(1,B)\).

Begin with \(\psi,B\ge0\). For regular compact support, Lemma~\ref{latticemarginal} supplies the joint limit with an arbitrary compactly supported lattice test. We use it to localize at \(L_{\rm std}\), then control the singular part.

\stage{Specializing to the free lattice}
Choose a relatively compact identity neighborhood \(\mathcal U\subset G\) mapping injectively to \(\mathcal X\), symmetric under inversion and conjugation by every required reflection. Let \(\rho\ge0\) be a continuous bump on its image, \(\mu(\rho)=1\). Exterior powers of a frame and its inverse are uniformly bounded on the patch. The integral exterior lattice is discrete in every rank, hence every nonzero primitive submodule has covolume bounded below there. Choose \(c\) below this bound. For \(h\in\mathcal U\) conjugation is then a bijection from the cyclic matrices on \(L_{\rm std}=\Oa^n\) to those on \(hL_{\rm std}\), all passing the cuts. Set
\[
 B^-_{\mathcal U}(Z)=\inf_{h\in\overline{\mathcal U}} B(hZh^{-1}),\qquad
 B^+_{\mathcal U}(Z)=\sup_{h\in\overline{\mathcal U}} B(hZh^{-1}).
\]
For a small patch these envelopes are continuous on a common regular compact support and have the required symmetries. On that support all contributing \(x_T\) lie in a fixed compact set \(K_0\), enlarged slightly if necessary. In fact on the patch
\[
 |\log H_j(hL_{\rm std})-\log H_j(L_{\rm std})|\le C_{\mathcal U},
 \qquad \|x_T(hL_{\rm std},hXh^{-1})-x_T(L_{\rm std},X)\|_\infty
        \le 2C_{\mathcal U}/\log T ,
\]
where \(H_j\) refers to the corresponding conjugate prefix. With \(\omega_{\psi,K_0}\) the modulus of continuity there, set \(\epsilon_T=\omega_{\psi,K_0}(2C_{\mathcal U}/\log T)\). The two envelope inequalities and \(\mu(\rho)=1\) give
\[
 P^f_{T,c}(\rho\psi B^-_{\mathcal U})-e_T
       \le M_T(\psi,B)\le P^f_{T,c}(\rho\psi B^+_{\mathcal U})+e_T,
 \qquad 0\le e_T=\epsilon_T P^f_{T,c}(\rho B^+_{\mathcal U})=o(1).
\]
The unweighted regular limit bounds the last integral. First let \(T\to\infty\) and then shrink the patch. The envelopes converge uniformly to \(B\), so the regular-support assertion follows for every bounded continuous weight.

\stage{Controlling the singular part at finite height}
Write \(\Sigma_u=\{Y:\min_v\|Y_v^{\,n-1}\|_{\Frob}\le u\}\). On a fixed identity patch conjugation enlarges the box and \(u\) by bounded factors. For a prefix of \(K\)-dimension \(r_j\), put \(k=qr_j\) and \(N=qn\). Choose \(N-k\) independent rows \(a_i\) of the integral matrix representing \(g_j(X)\) before changing the frame. Its saturated kernel has primitive normal wedge obtained by dividing \(a_1\wedge\cdots\wedge a_{N-k}\) by the gcd of its minors. With uniform metric distortion,
\[
 1\ll H_j\ll\|a_1\wedge\cdots\wedge a_{N-k}\|
       \ll_M(1+T)^{r_j(N-k)} .
\]
Thus even on singular support \(x_T\) lies in a fixed bounded range \(\Omega_1\). Bound the patch bump by its supremum and drop the ideal cuts while retaining \(\Omega_1\).

We now track the change of the power \(Z^{n-1}\) in the rough domination of Lemma~\ref{joint}. At one place in an adapted orthonormal frame write \(B_j=Z_{jj}\), \(C_{jl}=Z_{jl}\). The eigenvalues of \(B_j\) are the fixed layer eigenvalues divided by \(T\). On a bounded ball the characteristic equation gives, uniformly in the frame,
\[
 B_j^{b_j}=-\sum_{a=1}^{b_j}c_{ja}(T)B_j^{b_j-a},\qquad
 |c_{ja}(T)|\ll_{K,f}T^{-a},\qquad \|B_j^{b_j}\|_{\Frob}\ll_{M,K,f}T^{-1}.
\]
In expanding \(Z^{n-1}\), a word visiting distinct blocks \(i_1<\cdots<i_p\) and avoiding every \(B_j^{b_j}\) has length at most
\((p-1)+\sum_{a\le p}(b_{i_a}-1)=\sum_{a\le p}b_{i_a}-1\).
Equality at \(n-1\) forces the word
\[
 B_1^{b_1-1}C_{12}B_2^{b_2-1}\cdots C_{\ell-1,\ell}B_\ell^{b_\ell-1}
\]
at entry \(1,\ell\), also valid for \(\ell=1\). Every other term is \(O_M(T^{-1})\) for fixed \(K,f\). Put \(s_*=\prod_{e\in\mathcal E}s_e\). Scaling the exceptional blocks by \(1/2\le s_e\le2\) gives
\begin{equation}\label{finitesingular}
 \|(\mathcal R_s Z)^{n-1}-s_* Z^{n-1}\|_{\Frob}\ll_{M,K,f}T^{-1}.
\end{equation}
For example, \(Z=T^{-1}\left(\begin{smallmatrix}1&T\\0&1\end{smallmatrix}\right)\) with \(f=(t-1)^2\) gives \(\mathcal R_s Z-s Z=(1-s)I_2/T\); exact homogeneity already fails in this two-dimensional fiber.

In the proof of the rough domination use the actual flag frame for each smoothed block. Its diagonal spectra remain prescribed throughout. By \eqref{finitesingular} every image of a contributing original matrix lies in a fixed enlarged ball and in \(\Sigma_{C(u+T^{-1})}\). For each fixed \(u>0\) and sufficiently large \(T\) one may therefore use the closed set \(\{\|Y\|\le M'\}\cap\Sigma_{C'u}\) as the outer test in that proof. Formula \eqref{continuousfactor}, bounded \(|\Omega_1|\), and the upper bound for a compact closed set under vague convergence give
\[
 \limsup_{T\to\infty}M_T(\mathbf1_{\{H_\Box\le M\}\cap\Sigma_u})
       \ll_{M,K,f}\eta_n(\{\|Y\|\le M'\}\cap\Sigma_{C'u}).
\]
To evaluate the last bound, write \(Y_v=kUk^{-1}\) in the angular integral. Then \(\|Y_v^{n-1}\|_{\Frob}=\prod_{a<n}|U_{a,a+1}|\). Bound every upper entry by \(R\ge1\). With \(q_v=[K_v:\R]\), radial coordinates \(z_a=(|U_{a,a+1}|/R)^{q_v}\) on these disks give constant times \(dz_a\). For \(p=n-1\), the substitution \(z_j=e^{-t_j}\) yields, for \(0<a\le1\),
\[
 \vol\{z\in[0,1]^p:\textstyle\prod_j z_j\le a\}
 =\frac1{(p-1)!}\int_{\log(1/a)}^\infty e^{-t}t^{p-1}\,dt
 =a\sum_{j=0}^{p-1}\frac{(\log(1/a))^j}{j!}.
\]
Take \(a=(C'u/R^p)^{q_v}\), enlarging \(R\) if necessary. The other angular entries are bounded. Summing over the offending places gives
\begin{equation}\label{sectordecay}
 \limsup_{T\to\infty}M_T(\mathbf1_{\{H_\Box\le M\}\cap\Sigma_u})
 \ll_{M,K,f}\sum_{v\mid\infty} u^{q_v}(1+|\log u|)^{n-2}
 \qquad(0<u\le1/2).
\end{equation}
The constant is independent of \(u\); \(u\) is fixed before taking the upper limit in \(T\).

\stage{Removing the cutoff}
Let \(\chi_u(Y)\in[0,1]\) be continuous in \(\min_v\|Y_v^{n-1}\|_{\Frob}\), zero up to \(u\) and one from \(2u\) onwards. Its product with \(B\) has regular compact support and preserves the required reflections. For \(\operatorname{supp}B\subseteq\{H_\Box\le M\}\),
\[
 |M_T(\psi,B)-M_T(\psi,\chi_u B)|
 \le\|\psi\|_\infty\|B\|_\infty M_T(\mathbf1_{\{H_\Box\le M\}\cap\Sigma_{2u}}).
\]
Use the regular-support formula, \eqref{sectordecay}, and the same angular bound for \(\eta_n\), then let \(u\downarrow0\). Positive and negative parts give the assertion for signed tests.
Finally let \(\Omega\) be bounded with \(\eta_n(\partial \Omega)=0\). For \(\varepsilon>0\) regularity supplies compactly supported continuous functions \(0\le B_-\le\mathbf1_{\Omega}\le B_+\le1\) with \(\eta_n(B_+-B_-)<\varepsilon\). When Theorem~\ref{law} requires symmetry, average both envelopes over the finite group of real-place reflections. Invariance of \(\Omega\) preserves the inequalities and invariance of \(\eta_n\) preserves the error. In particular
\[
 \limsup_{T\to\infty}|M_T(\psi,\mathbf1_{\Omega})-M_T(\psi,B_-)|
 \le\|\psi\|_\infty A_0\vol(P)\eta_n(B_+-B_-)
 <\|\psi\|_\infty A_0\vol(P)\varepsilon .
\]
The limiting angular integrals obey the same estimate. Let \(\varepsilon\downarrow0\).
\end{proof}

\noindent To identify \(A_0\), apply Proposition~\ref{localized} to \(g\) and the product of local Frobenius balls with \(\psi=1\). Every matrix on this squarefree fiber is cyclic. Cai--Deng's coefficient gives \(A_0\vol(Z)\prod_v\kappa_{q_vD}=C^\Frob_{K,g}\), so \(A_0=a_{K,g}\). This proves Theorem~\ref{law} before any estimate for noncyclic matrices is used.

\begin{example}[The mixed-sign proportion for \(t^2+1\)]
Keep \(K=\Q(\sqrt3)\), \(s=\sqrt3\), \(F=K(i)\), \(\alpha=(s+i)/2\), and the real places \(v_+,v_-\). For \(f=t^2+1\) the mixed-sign proportion is
\[
 \lim_{T\to\infty}
 \frac{\#\{X\in\mathcal X_f(\Oa_K):H_\Box(X)\le T,\
              (X_{v_+})_{21}(X_{v_-})_{21}<0\}}
      {N_f^\Box(T)}=\frac3{10}.
\]

\stage{The ten normalized lattices}
Put \(O=\Oa_F=\Oa_K[\alpha]\) and \(R=\Oa_K[i]=\Oa_K+2O\).
The preceding example proves \(h_F=1\), and \(\sqrt{12}/2<2\) gives \(h_K=1\).
Integral conjugacy classes correspond to full fractional \(R\)-ideals modulo \(F^\times\).
Every such ideal is free over \(\Oa_K\), and
\[
 OI=O\quad\Longrightarrow\quad 2O\subseteq I\subseteq O,\qquad
 I\sim I'\ \Longleftrightarrow\ I'=uI\quad(u\in O^\times)
\]
for normalized ideals. Let \(\pi=1+\sqrt3\), \(\mathfrak p=(\pi)\), and \(\Oa_2=\Oa_K/\mathfrak p^2\). Then
\[
 N\pi=-2,\quad \pi^2=2(2+\sqrt3),\quad
 (2)=\mathfrak p^2,\quad \Oa_K/\mathfrak p=\mathbb F_2,\quad
 O/\mathfrak pO=\mathbb F_4 .
\]
Thus \(I/2O\subset O/2O\simeq\Oa_2^2\) must contain a primitive vector.
Such a vector is a unit modulo \(2O\), so this condition is equivalent to \(OI=O\).
Every preimage is \(R\)-stable. The elementary divisors and their counts are
\[
\begin{array}{c|c|c|c}
 j&\text{elementary divisors of }I\subset O&
        \text{choice of }I/2O&a_j\\ \hline
 0&(1,1)&\Oa_2^2&1\\
 1&(1,\mathfrak p)&\text{preimage of a line in }\mathbb F_2^2&
             (2^2-1)/(2-1)=3\\
 2&(1,\mathfrak p^2)&\text{free rank-one summand of }\Oa_2^2&
             (4^2-2^2)/2=6
\end{array}
\]

\stage{The sign is determined by \(j\)}
In the basis \(1,\alpha\), multiplication by \(i\) is \(M=2C-sI_2\), where
\(C=\left(\begin{smallmatrix}0&-1\\1&s\end{smallmatrix}\right)\).
If the columns of \(h\) form an \(\Oa_K\)-basis of \(I\), then
\[
 X_I=h^{-1}Mh,\qquad \det h=\epsilon\pi^j,\quad\epsilon\in\Oa_K^\times.
\]
For \(z=(u,w)^t\ne0\) at either real place,
\[
 \det(z,M_vz)=2\big((u+s_vw/2)^2+w^2/4\big)>0 .
\]
Using \(z=h_ve_1\) and \(N\epsilon=1\) therefore gives
\[
 \operatorname{sgn}\big((X_{I,v_+})_{21}(X_{I,v_-})_{21}\big)
 =\operatorname{sgn}N(\det h)=(-1)^j.
\]
This holds throughout the integral conjugacy class.

\stage{Orbit weights and localization}
Set \(U_I=\{u\in O^\times:uI=I\}\) and \(w_I=[O^\times:U_I]\).
Write \(P_T^{[I]}\) for the class contribution to the averaged measure.
Since \(\ell=1\), unfolding and the finite covering \(G/U_I\to G/O^\times\) give
\[
 P_T^{[I]}(B)=w_I P_T^{[O]}(B),\qquad
 \sum_{\substack{[I]\\\text{index level }j}}w_I=a_j,\qquad
 \sum_{[I]}w_I=10.
\]
Here an archimedean change of basis identifies the orbital integrals, and
\(w_I\) is also the size of the \(O^\times\)-orbit of a normalized lattice.

For a reflection-invariant regular test \(B\), the same identity gives the
class marginal \(w_I A_0\eta_2(B)/10\).
Centralizer translations, scalars and reflections preserve the integral
\(R\)-module class. Its short-vector and singular-set errors are bounded by
the corresponding total positive errors. Lemmas~\ref{flowinvariance}--\ref{latticemarginal}
and the class-preserving identity patch in Proposition~\ref{localized} give
\[
 \begin{aligned}
 \lim_{T\to\infty}P_T^{[I]}(\rho B)
    &=\frac{w_I}{10}\mu(\rho)A_0\eta_2(B),\\
 \lim_{T\to\infty}\frac{N_{[I]}^\Box(T)}{N_f^\Box(T)}
    &=\frac{w_I}{10},\\
 \lim_{T\to\infty}\frac{N_{\rm mixed}^\Box(T)}{N_f^\Box(T)}
    &=\frac{\sum_{[I]:\,j=1}w_I}{\sum_{[I]}w_I}
      =\frac{3}{1+3+6}.
 \end{aligned}
\]
The second line uses the singular-set bound and the null boundary of the entry box.
The sign condition enters through the class calculation.

For comparison, a real reflection exchanges the two sign sets and preserves
the box under \(\eta_2\). Their common boundary is null, since
\((k(uE_{12})k^{-1})_{21}=0\) only for \(u=0\) or \(ke_1\parallel e_1\).
Thus the angular measure assigns proportion \(1/2\) to mixed signs.
\end{example}

\section{Rational conjugacy types and the unrestricted count}\label{types}
The cyclic calculation counts the top rational type. We now bound the other types. The new issue is that two successive layers for the same irreducible polynomial can have different multiplicities. Their intertwiner is then a rectangular matrix over a field.

\subsection{A mean estimate for a Gaussian on rectangular matrices}
Let \(E\) be a number field, write \(m=[E:\Q]\), and put
\(H_a=\{h\in\GL_a(E_\infty):|N\det h|=1\}\).
\begin{lem}[Full-column-rank mean]\label{gaussmean}
Let \(1\le b\le a\), and let \(I\) be a full lattice in \(M_{a,b}(E)\). Suppose \(\Delta\subset H_a\) preserves \(I\) under left multiplication and is commensurable with \(\GL_a(\Oa_E)\) in the chosen coordinates. Give \(H_a/\Delta\) invariant probability. For every positive definite real quadratic form \(Q_0\) on \(M_{a,b}(E_\infty)\),
\[
 \int_{H_a/\Delta}\ \sum_{\substack{z\in I\\\operatorname{rank}_Ez=b}}
       e^{-Q_0(hz)}\,dh
 \ \le C(E,I,\Delta)\int_{M_{a,b}(E_\infty)} e^{-Q_0(X)}\,dX .
\]
The constant is independent of the eccentricity and scale of \(Q_0\).
\end{lem}
\begin{proof}
Put \(\Gamma_a=\GL_a(\Oa_E)\) and choose an integer \(d\ge1\) with
\(dI\subseteq M_{a,b}(\Oa_E)\).
Pass to the common finite cover \(H_a/(\Delta\cap\Gamma_a)\), always using probability.
The pointwise inclusion gives
\[
 \sum_{\substack{z\in I\\\operatorname{rank}z=b}}e^{-Q_0(hz)}
 \le\sum_{\substack{z\in M_{a,b}(\Oa_E)\\\operatorname{rank}z=b}}
                  e^{-Q_0(hz/d)} .
\]
The right-hand sum descends to \(H_a/\Gamma_a\), and
\[
 \int_{M_{a,b}(E_\infty)}e^{-Q_0(X/d)}\,dX
       =d^{mab}\int_{M_{a,b}(E_\infty)}e^{-Q_0(X)}\,dX.
\]
It suffices to prove the assertion for the standard model.

\stage{The case \(b<a\)}
Let \(\mathbb A_E\) be the adeles and \(\mathcal X_a=\GL_a(\mathbb A_E)^1/\GL_a(E)\) with its probability measure. For \(b<a\) the independent-column formula \cite[(1.5)]{AS}, \cite[Theorem 1.2]{Kim} is
\[
 \int_{\mathcal X_a}\sum_{\substack{z\in M_{a,b}(E)\\\operatorname{rank}_Ez=b}}
        \Phi(hz)\,dh
    =\int_{M_{a,b}(\mathbb A_E)}\Phi(X)\,dX ,
\]
with adelic additive measure normalized by \(\vol(\mathbb A_E/E)=1\). Take
\(\Phi=e^{-Q_0}\mathbf1_{M_{a,b}(\widehat{\Oa}_E)}\) as an archimedean times finite test. Write \(p_0>0\) for the mass of the component represented by \(H_a\times\GL_a(\widehat{\Oa}_E)\). On it the finite test selects the integral columns and the conditional archimedean probability is the required lattice probability. Since \(\Phi\ge0\), the standard lattice mean \(\mathcal M(Q_0)\) satisfies
\[
 p_0\mathcal M(Q_0)\le \int_{M_{a,b}(\mathbb A_E)}\Phi\,dX
     = c_{E,a,b}\int_{M_{a,b}(E_\infty)}e^{-Q_0(X)}\,dX .
\]
Both constants are independent of the shape of \(Q_0\).

\stage{The case \(b=a\)}
Two nonsingular integral square matrices lie in the same left \(\GL_a(\Oa_E)\)-orbit precisely when they have the same row lattice. Its ordered rows are a basis of that lattice, so the stabilizer is trivial. For a representative \(z\), the index of its row lattice is \(t=|N\det z|\). Right translation on \(H_a\) shows that the unfolded orbital integral has the form
\[
 F(t)=\int_{H_a} e^{-Q_0(t^{1/(am)}h)}\,dh.
\]
Indeed \(t^{-1/(am)}z\in H_a\). The function \(F\) is decreasing since \(Q_0(t^{1/(am)}h)=t^{2/(am)}Q_0(h)\).

The number \(A(y)\) of such row lattices of index at most \(y\) is \(O_{E,a}(y^a)\). To see this, majorize by all full-index submodules; successive coordinate-ideal projection gives their index series \(\prod_{j=0}^{a-1}\zeta_E(s-j)\). Equivalently at each finite prime this is the triangular-basis calculation over the local discrete valuation ring. Writing \(c(m)\) for the coefficients of the first \(a-1\) factors, ideal counting bounds the last convolution by
\[
 A(y)\le\sum_{m\le y}c(m)\sum_{N\mathfrak a\le y/m}(N\mathfrak a)^{a-1}
 \ll y^a\sum_{m\ge1}\frac{c(m)}{m^a}
 =y^a\,O_{E,a}\!\left(\prod_{j=0}^{a-2}\zeta_E(a-j)\right).
\]
The empty product for \(a=1\) is one. Write \(r(u)=\sup\{t>0:F(t)>u\}\), with \(\sup\varnothing=0\). Integrating superlevel sets gives
\[
 \int F(t)\,dA(t)\le\int_0^\infty A(r(u))\,du
 \ll\int_0^\infty r(u)^a\,du
 =a\int_0^\infty F(t)t^{a-1}\,dt .
\]
Haar measure on \(\GL_a(E_\infty)\) is a constant times \(dh\,dt/t\) in these coordinates. Lebesgue measure on its matrix entries is \(t^a\) times this Haar measure up to a fixed normalization, since left multiplication has real determinant \(|N\det|^a\). The singular matrices have entry measure zero. The last integral is thus a fixed multiple of the Gaussian entry integral for every \(Q_0\). This proves the square case.
\end{proof}

\subsection{The bound for a rational type}
For \(X\in\mathcal X_f(\Oa)\) put \(K_{i,0}=0,\ K_{i,j}=\ker_K p_i(X)^j\) and
\[
 \begin{aligned}
 r_{ij}&=\frac{\dim_K K_{i,j}-\dim_K K_{i,j-1}}{d_i},
       &s&=\#\{(i,j):j\ge1,\ r_{ij}>0\}.
 \end{aligned}
\]
The type is \(\operatorname{type}(X)=(r_{ij})_{1\le i\le k,\,j\ge1}\), and its counting function is
\[
 N_{\boldsymbol r}(T)=\#\{X\in\mathcal X_f(\Oa):\operatorname{type}(X)=\boldsymbol r,\
 H_\Box(X)\le T\}.
\]
Equivalently, if the \(p_i\)-primary part decomposes as \(\bigoplus_a K[t]/p_i^{\lambda_{ia}}\), then
\[
       r_{ij}=\#\{a:\lambda_{ia}\ge j\},\qquad
       r_{i1}\ge r_{i2}\ge\cdots,\qquad \sum_{j\ge1}r_{ij}=e_i .
\]
\begin{lem}[The centralizer dimension of a rational type]\label{typedimension}
For a matrix of type \(\boldsymbol r\),
\[
 \dim_K\End_{K[t]}(K^n)=\sum_{i,j}d_i r_{ij}^2,\qquad
 \dim(\GL_n\cdot X)=n^2-\sum_{i,j}d_i r_{ij}^2 .
\]
\end{lem}
\begin{proof}
Different primaries have no nonzero intertwiners by B\'ezout's identity. Within one primary, an intertwiner is determined by the image of \(1\), so
\[
 \Hom_{K[t]}(K[t]/p_i^a,K[t]/p_i^b)
 \simeq\{z\in K[t]/p_i^b:p_i^az=0\},\qquad
 \dim_K\Hom=d_i\min(a,b).
\]
For \(M_i=\bigoplus_h K[t]/p_i^{\lambda_{ih}}\), sum the dimensions:
\[
 \dim_K\End_{K[t]}M_i
  =d_i\sum_{h,k}\min(\lambda_{ih},\lambda_{ik})
  =d_i\sum_{j\ge1}\left(\sum_h\mathbf1_{\lambda_{ih}\ge j}\right)^2
  =d_i\sum_jr_{ij}^2 .
\]
The invertible centralizers form a nonempty Zariski-open subset of the endomorphism algebra, with the same dimension. The orbit--stabilizer dimension formula gives the second assertion.
\end{proof}

\begin{prop}[Counting a fixed rational type]\label{strata}
Fix $K,f$. For the type with column lengths $r_{ij}$ and $s$ nonzero columns in total,
\begin{equation}\label{type}
 N_{\boldsymbol r}(T)\ll_{K,f}
 T^{q(n^2-\sum_{i,j}d_i r_{ij}^2)/2}(\log(2T))^{s-1}.
\end{equation}
\end{prop}
\begin{proof}
\stage{The canonical filtration}
Order the primaries as in Section~\ref{flags}. If \(V_i=\ker p_i(X)^{e_i}\), the successive terms of the canonical filtration are
\[
 F_{i,j}=\Big(\bigoplus_{h<i}V_h\Big)\oplus K_{i,j},\qquad
 F_{i,j}/F_{i,j-1}\simeq E_i^{\,r_{ij}},\qquad E_i=K[t]/p_i ,
\]
where each nonzero layer appears once. The operator on the quotient is multiplication by \(t\bmod p_i\). Integrate the fixed-lattice count over the identity patch with a probability bump and enlarge the box by its bounded conjugation distortion. For
\(q_{i,j}(t)=(\prod_{h<i}p_h(t)^{e_h})p_i(t)^j\),
the kernel-minor estimate of Proposition~\ref{localized} gives
\[
 F_{i,j}=\ker q_{i,j}(X),\qquad
 1\ll h_{F_{i,j}}\ll(1+T)^{C_{K,f}},\qquad
 |\log h_{F_{i,j}}|\ll\log(2T).
\]
Bounding the bump by its supremum and integrating the \(s-1\) free heights therefore costs
\[
 \int_{\text{retained heights}}\prod_{a<s}d\log H_a
       \ll_{K,f}(\log(2T))^{s-1}.
\]

For one layer write \(R=\Oa[t]/p_i\), \(E=K[t]/p_i\), \(a=r_{ij}\). For each full \(R\)-lattice \(I\subset E^a\) and a fixed nonzero integer \(d\) with \(d\Oa_E\subset R\),
\[
 d(\Oa_E I)\subseteq I\subseteq \Oa_E I .
\]
Up to \(\GL_a(E)\), the saturations \(\Oa_E I\) range over finitely many Steinitz classes, and each quotient \(\Oa_E I/d(\Oa_E I)\) is finite. This gives finitely many models and their arithmetic stabilizers are commensurable with \(\GL_a(\Oa_E)\).
\stage{Shears and the remaining full-rank entries}
Fix a run for \(p_i\), put \(E=K[t]/p_i\), and let \(s_i=\#\{j:r_{ij}>0\}\). Its layers are \(V_j=E^{m_j}\) with \(m_j=r_{ij}\) and \(m_1\ge\cdots\ge m_{s_i}\). The positive real scale on layer \(j\) is still \(a_j\) from \eqref{jac}. If \(K_j=\ker p_i(X)^j\) in this primary, then
\[
 p_i(X):K_{j+1}/K_j\hookrightarrow K_j/K_{j-1},
 \qquad \operatorname{rank}_E z_j=m_{j+1},
\]
where \(p_i'(A_j)z_j\) is the induced map after identifying the semisimple layers. Conversely, on a filtration \(0=W_0\subset\cdots\subset W_k\) let \(N=p_i(X)\) lower layers and write \(\overline N_h:\operatorname{gr}_h W\to\operatorname{gr}_{h-1}W\) for the adjacent injections. For \(h>j\),
\[
 \operatorname{gr}_h(N^j)=\overline N_{h-j+1}\cdots\overline N_h
 \text{ is injective},\qquad
 w\in W_h\setminus W_{h-1}\ \Longrightarrow\ N^j w\ne0.
\]
Since \(N^jW_j=0\), this proves \(\ker N^j=W_j\) and identifies the canonical filtration uniquely. At distance \(h\ge2\) the projected offset map has components \(L_h(v)_j=v_jz_{j+h-1}-z_jv_{j+1}\). Injectivity gives an \(E\)-linear left inverse \(\ell_k:V_k\to V_{k+1}\), \(\ell_k z_k=1\), for each \(z_k\). Given \(y_j:V_{j+h}\to V_j\), set the last source variable to zero and solve from right to left:
\[
 v_j=(y_j+z_jv_{j+1})\ell_{j+h-1}
 \quad\Longrightarrow\quad v_jz_{j+h-1}-z_jv_{j+1}=y_j .
\]
This is a right inverse on vector spaces, hence also after extension to the archimedean factors. The left inverses need no integrality. Separately, composition by an integral intertwiner preserves the projected periods, as in Lemma~\ref{sheardist}. Thus \(L_h\) descends to a surjection of the period tori and pushes probability Haar to probability Haar. Descending induction in \(h\) leaves just the full-rank conditions on adjacent blocks.

\stage{Integration along a sequence of equal factors}
In an orthonormal flag frame let \(X_{jl,v}\) be the physical blocks. On \(H_\Box(X)\le T\) one has \(\|X_v\|_{\Frob}\le nT\) and \(\sum_{j,l}\|X_{jl,v}\|_{\Frob}^2=\|X_v\|_{\Frob}^2\). Consequently
\[
 \mathbf1_{H_\Box(X)\le T}\le e^{n^2|S_\infty|}
   \prod_{j,v}\mathbf1_{\{\|X_{jj,v}\|_{\Frob}\le nT\}}\,
   \prod_{j<l}\exp\!\left(-T^{-2}\sum_v\|X_{jl,v}\|_{\Frob}^2\right).
\]
The diagonal cutoffs will bound the orbital integrals after the upper blocks have been integrated. Fix the right centralizer and both orbital coordinates of one adjacent equal-primary edge. In the untransformed block the test is \(e^{-Q(z+w)}\) for a positive definite form \(Q\), where \(w\in W\) and \(z\) is the intertwiner. In compatible Lebesgue coordinates, complete the square:
\[
 \begin{aligned}
 Q(z+w)&=z^tQ_Vz+2z^tQ_{VW}w+w^tQ_Ww,\\
 Q_0(z)&=z^t(Q_V-Q_{VW}Q_W^{-1}Q_{VW}^t)z>0\quad(z\ne0),\\
 Q(z+w)&=Q_0(z)+
 (w+Q_W^{-1}Q_{VW}^tz)^tQ_W(w+Q_W^{-1}Q_{VW}^tz).
 \end{aligned}
\]
With \(w_0=\dim_\R W\), integration gives
\[
 \int_W e^{-Q(z+w)}\,dw
 =c(Q,W)e^{-Q_0(z)},\qquad
 c(Q,W)=\pi^{w_0/2}(\det Q_W)^{-1/2}.
\]
When \(W=0\), take \(Q_0=Q\) and \(c(Q,W)=1\). The left centralizer \(H_{m_j}\) on an edge \(j,j+1\) acts on \(V=M_{m_j,m_{j+1}}(E_\infty)\) by \(z\mapsto hz\). Each \(h\in H_{m_j}\) preserves \(W\) and has real determinant of modulus one on \(V\) and on the entire block, hence on \(W\). Changing \(w\) absorbs its action on \(W\). Let \(I\) be the projected lattice, \(\Delta\) its left arithmetic centralizer, \(a=m_j\) and \(b=m_{j+1}\). Writing \(h_W=\covol(\Lambda\cap W)\), the edge integral is bounded by
\[
 \begin{aligned}
 \frac{c(Q,W)}{h_W}
 \int_{H_a/\Delta}\sum_{\substack{z\in I\\\operatorname{rank}_E z=b}}
       e^{-Q_0(hz)}\,dh
 &\le \frac{C(E,I,\Delta)c(Q,W)}{h_W}
               \int_V e^{-Q_0(z)}\,dz\\
 &= C(E,I,\Delta)\covol(I)
          \int_{V\oplus W} e^{-Q(C')}\frac{dC'}{\covol\Lambda}.
 \end{aligned}
\]
The inequality is Lemma~\ref{gaussmean}; \(\covol\Lambda=h_W\covol I\) gives the equality. The projected lattice is stable under \(\Delta\), commensurable in an \(E\)-basis with \(\GL_a(\Oa_E)\). For \(b<a\) the lemma uses the independent-column mean, and for \(b=a\) it unfolds row lattices. Only the fixed finite models enter the multiplier, even as the right shape and \(Q_0\) vary.

Integrate first the ordinary blocks, whose determinants remove both centralizer variables. In each run integrate the remaining edges from left to right. At the first edge the left centralizer occurs nowhere else; replacing that edge by ordinary entry measure removes its right centralizer from that edge as well. Induction applies to the next pair. The uniformity in Lemma~\ref{gaussmean} is what permits this order of integration. Indexing all layers in order, let \(b_j\) be the \(K\)-dimension of layer \(j\), its primary degree times its column multiplicity. The scale factors thus cancel:
\[
 \underbrace{\prod_{j<l}(a_j/a_l)^{qb_jb_l}}_{\text{flag modulus}}\,
 \underbrace{\prod_{j<l}(Ta_l/a_j)^{qb_jb_l}}_{\text{Gaussian entry integrals}}
     = T^{q\sum_{j<l}b_jb_l}.
\]

\stage{The exponent}
A semisimple layer of size $b=da$ contributes orbital volume $O(T^{q(b^2-da^2)/2})$. At a complex place the spectral blocks are scalar and the rectangular entries count twice in real dimension.

At a real place a conjugate pair \(\alpha\pm i\beta\) of multiplicity \(a\) contributes the symmetric space \(\GL_{2a}(\R)/\GL_a(\C)\). Real positive central scalars are in its stabilizer. On the positive determinant component normalize to \(\SL_{2a}(\R)\) with stabilizer \(\{h\in\GL_a(\C):|\det_\C h|=1\}\); the other component has the same bound. Write the polar parameter as \(\operatorname{diag}(t_1,\ldots,t_a,-t_1,\ldots,-t_a)\) with \(t_1\ge\cdots\ge t_a\ge0\). Conjugating \(J=\left(\begin{smallmatrix}0&-I\\I&0\end{smallmatrix}\right)\) by its exponential changes the upper right block to \(-\operatorname{diag}(e^{2t_i})\) and the lower left to \(\operatorname{diag}(e^{-2t_i})\). Therefore the height cutoff on \(\alpha I+\beta J\) gives \(e^{2t_1}\ll_{\beta} T\), where the enclosing ball is fixed after scaling by \(T\).

For clarity, in the polar density bound \cite[(3.7)--(3.8)]{GOS} each positive restricted root \(\lambda\) contributes at most \(O(e^{m_\lambda\lambda(t)})\). The differences of the \(2a\) diagonal weights give \(t_i-t_j,t_i+t_j\) twice for \(i<j\), and \(2t_i\) once. Consequently
\[
 \begin{aligned}
 \sum_{\lambda>0}m_\lambda\lambda(t)
 &=2\sum_{i<j}\big((t_i-t_j)+(t_i+t_j)\big)+2\sum_i t_i\\
 &=2\sum_{i=1}^a(2a+1-2i)t_i .
 \end{aligned}
\]
Enlarging the ordered chamber to the enclosing cube gives
\[
 \int_{\mathrm{height}\le T}d\mu_{\rm orb}
 \ll\prod_{i=1}^a\int_0^{\log(CT)/2}e^{2(2a+1-2i)t}\,dt
 \ll T^{\sum_{i=1}^a(2a+1-2i)}=T^{a^2}.
\]
For \(u\) real roots and \(v\) conjugate pairs, the spectral block widths \(m_\nu\) are \(a\) and \(2a\), respectively. The real exponent is therefore
\[
 \sum_{\nu<\mu}m_\nu m_\mu+va^2
 =\tfrac12\!\left(b^2-(u+4v)a^2\right)+va^2
 =\tfrac12(b^2-da^2),\qquad d=u+2v .
\]
At a complex place there are \(d\) widths \(a\), with two real coordinates per complex entry. This gives the local exponent \(q_v(b^2-da^2)/2\) at every place.

The entry integrals, the layer volumes and the \(s-1\) height integrals now give
\[
 \begin{aligned}
 N_{\boldsymbol r}(T)
 &\ll T^{q\sum_{j<k}b_jb_k}
       T^{\frac q2(\sum_j b_j^2-\sum_{i,j}d_i r_{ij}^2)}
       (\log(2T))^{s-1}\\
 &=T^{\frac q2(n^2-\sum_{i,j}d_i r_{ij}^2)}
       (\log(2T))^{s-1}.
 \end{aligned}
\]
Here \(2\sum_{j<k}b_jb_k+\sum_jb_j^2=(\sum_jb_j)^2=n^2\).
The exponent deficit is
\[
 qD_n-\frac q2\left(n^2-\sum_{i,j}d_i r_{ij}^2\right)
     =q\sum_{i,j}d_i\binom{r_{ij}}2>0
\]
for every noncyclic type.
\end{proof}

\begin{proof}[Proof of Corollary~\ref{small}]
A noncyclic type has some \(r_{ij}\ge2\) with \(e_i\ge2\). Hence
\[
 qD_n-\tfrac q2\big(n^2-\sum_{i,j}d_i r_{ij}^2\big)
   =q\sum_{i,j}d_i\binom{r_{ij}}2\ge q\delta,\qquad
 s\le\sum_{i,j}r_{ij}-1=\ell-1.
\]
Substitute in \eqref{type} and sum over the finitely many primary partitions.
\end{proof}

\begin{proof}[Proof of Theorem~\ref{main}]
For a fixed separation \(g\), Proposition~\ref{localized} with \(\psi=1\), together with the negligible noncyclic count, gives
\[
 \begin{aligned}
 \lim_{T\to\infty}\frac{N_f^\Box(T)}{S(T)}
   &=A_0\vol(P)\prod_vJ_{n,v},\\
 \lim_{T\to\infty}\frac{N_g^\Box(T)}{S(T)}
   &=A_0\vol(Z)\prod_vJ_{n,v},\\
 C^\Frob_{K,g}
   &=A_0\vol(Z)\prod_v\kappa_{q_vD}.
 \end{aligned}
\]
The third line evaluates the same squarefree count as \cite[Theorem 1.2]{CD}. Dividing the first line by the second yields \(\vol(P)/\vol(Z)=1/\prod_i e_i!\) by \eqref{ch}. Substitution of \(A_0\) from the third line gives \eqref{absolute}.
\end{proof}

\section{Arithmetic normalization}\label{arithsec}
\subsection{The number-field coefficient}
We express the Frobenius coefficient of the separation in the notation of \cite[Theorem 1.2]{CD}. Put \(A_g=K[t]/(g)=\prod_j E_j\), \(\alpha=t\bmod g\), and let \(r_2\) count the complex places of \(K\).
For an order \(R\) in a finite \'etale \(\Q\)-algebra, let \(R^\vee=\Hom_{\Z}(R,\Z)\), with \(R\)-action \((a\lambda)(r)=\lambda(ar)\), and define for \(\operatorname{Re}s\) large
\[
 \zeta_R(s)=\sum_{\substack{M\subset R^\vee,\ RM\subseteq M\\
                              [R^\vee:M]<\infty}} [R^\vee:M]^{-s},
\]
where \(M\) ranges over additive subgroups. Embed the factors \(E_j\) in a common splitting field as in \cite[Theorem 1.2]{CD}. Sum over common cyclic subfields \(K\subseteq E\subseteq \bigcap_j E_j\) unramified at every place, including infinity. With the resulting diagonal \(E\)-algebra structure on \(A_g\), put
\[
 \mathfrak S(g)=\sum_E\varphi([E:K])\,\delta(E)
    \lim_{s\to1}(s-1)^\ell\zeta_{\Oa_E[\alpha]}(s).
\]
For \(v\) finite write \(E_v=E\otimes_KK_v\), \(A_v=A_g\otimes_KK_v\), and use their maximal orders. Set
\[
 \Oa_{E_v}^1=\ker(N_{E_v/K_v}:\Oa_{E_v}^{\times}\to\Oa_{K_v}^{\times}),
 \qquad
 \delta(E)=\prod_{v\nmid\infty}
   [\Oa_{E_v}^1:\Oa_{E_v}^1\cap N_{A_v/E_v}(\Oa_{A_v}^{\times})].
\]
In particular \(\delta(K)=1\). Use the completion convention
\[
 \Lambda_K(s)=\zeta_K(s)\prod_{v\mid\infty}\Gamma_v(s),\quad
 \Gamma_{\R}(s)=\pi^{-s/2}\Gamma(s/2),\quad
 \Gamma_{\C}(s)=2(2\pi)^{-s}\Gamma(s).
\]
\begin{prop}[The arithmetic coefficient and the global ball]\label{arithmeticcoefficient}
With the preceding notation, the coefficient for the entry box is
\begin{equation}\label{arith}
 C^\Box_{K,f}=
 \frac{2^{r_2D} q^{\ell-1} n^{\ell-2}\prod_i d_i^{e_i}}
 {|\Delta_K|^{D}\pi^{nr_2}\operatorname*{Res}_{s=1}\Lambda_K(s)
                         \prod_{j=2}^n\Lambda_K(j)}
 \frac{\prod_v J_{n,v}}{\prod_i e_i!}\,\mathfrak S(g).
\end{equation}
Here \(\Delta_K\) is the field discriminant. For the single global Euclidean height \(H_{\rm eucl}(X)^2=\sum_v\|X_v\|_{\Frob}^2\), the coefficient is
\[
 C^{\rm eucl}_{K,f}=C^\Box_{K,f}\,
        \frac{\kappa_{qD}}{\prod_vJ_{n,v}} .
\]
\end{prop}
\begin{proof}
In \cite[Theorem 1.2]{CD}, the separated layer degrees are \(b_j=d_{i(j)}\), and the flag factor is \(q^{\ell-1}n^{\ell-2}\prod_jb_j\). The local convention \(2\,dx\,dy\) per complex upper entry contributes
\[
 w_n=2^{r_2D}\prod_{v\mid\infty}\kappa_{q_vD}.
\]
Substituting these factors and \(\mathfrak S(g)\) into the squarefree Frobenius coefficient, then applying \eqref{absolute}, cancels \(\prod_v\kappa_{q_vD}\) and gives \eqref{arith}.

In orthogonal or unitary triangular coordinates the sum of the squared Frobenius norms is the ordinary squared norm on \(\R^{qD}\). Hence
\[
 \eta_n(\{H_{\rm eucl}\le1\})=\kappa_{qD},\qquad
 \eta_n(\{H_\Box\le1\})=\prod_vJ_{n,v}.
\]
Both boundaries have zero angular mass by radial integration. Proposition~\ref{localized} supplies the cyclic coefficients. Since \(H_\Box\le H_{\rm eucl}\), Corollary~\ref{small} also bounds the noncyclic contribution in the global ball by \(o(S(T))\). If \(f\) is squarefree every matrix is cyclic, and \(n=1\) gives the single prescribed matrix directly. Taking the ratio of the angular masses proves the assertion.
\end{proof}

\subsection{The product over the rational simple orders}
\begin{proof}[Proof of Corollary~\ref{integercase}]
We evaluate the three factors in \eqref{marg}: the flag density, the diagonal density, and the volume of the ordered polytope.

\stage{Flag density} Over $\Z$ there is a single projective tuple for the degrees $b_1,\ldots,b_\ell$. Give the shape spaces their lattice probabilities, the frames invariant probability, and the shear quotient its invariant probability. For a primitive unoriented rank-$a$ sublattice of a unimodular lattice of rank $m=a+b$, Kim's mean value formula \cite{flagmean} gives the density
\[
 \frac{L_aL_b}{L_m}\,H^m\dd\log H
\]
in the sublattice covolume $H$. Indeed the coefficient is
\[
 \binom ma\prod_{i=1}^a
 \frac{\kappa_{m-i+1}\zeta^*(i)}
      {\kappa_i\zeta(m-i+1)},\qquad
 \zeta^*(1)=1,\quad \zeta^*(i)=\zeta(i)\ (i>1).
\]
Substitution of the ball volumes cancels the binomial coefficient and gives the displayed $L$-ratio. Conditional on the height, the normalized lattice and quotient shapes have invariant probabilities by parabolic disintegration. Put \(n_j=\sum_{a\ge j}b_a\). Iterating on the rescaled quotient telescopes:
\[
 \prod_{j<\ell}\frac{L_{b_j}L_{n_{j+1}}}{L_{n_j}}
      =\frac{\prod_{j=1}^\ell L_{b_j}}{L_n}.
\]
At each step the logarithmic height changes by a linear combination of earlier ones, with diagonal derivative one. Thus the height Jacobian is one, and the remaining height power is the modulus in \eqref{jac}. Every integral Hom lattice in these standard models has covolume one.

\stage{Density on one diagonal layer} For a single simple layer of degree $d$, Section~\ref{fixed} identifies its lattice-averaged orbital mass with the count on $\Z^d$ in the Frobenius ball. There is no proper ideal cutoff or logarithmic height coordinate in this case. By \cite[Theorem 1.1]{CD} that count has coefficient $\kappa_{D_d}A_p/L_d$. The trace-dual lattice of the monogenic order is $p'(\alpha)^{-1}\Z[\alpha]$, where $\alpha=t\bmod p$, so sublattices stable under the order give the ordinary ideal zeta function. Hence the summed diagonal multiplier of $\eta_d$ is $A_p/L_d$, also in degree one. The flag and diagonal constants together are $(\prod_i A_{p_i}^{e_i})/L_n$.

\stage{Volume of the chamber} To evaluate $\vol Z$, put $u_1,\ldots,u_\ell$ for the standard coordinate vectors. For $q=1$, the polytope in Section~\ref{flags} is the sum of segments
\[
   \sum_{j<k} b_jb_k[-\tfrac12,\tfrac12](u_j-u_k).
\]
Indeed the minimum of a subset sum on this zonotope is $-b(I)(n-b(I))/2$. For a linear functional with ordered coefficients, successive differences of coefficients express it as a multiple of the total sum plus a nonnegative combination of nested subset sums. Orienting each segment according to that order simultaneously attains the requisite extreme values. This proves equality of the support functions. In the first $\ell-1$ coordinates, a maximal minor of the edge vectors has absolute value one for a spanning tree and zero otherwise, by successively deleting leaves.

For \(\ell\ge2\), summing the maximal minors in the zonotope volume formula and then using the Pr\"ufer code gives
\[
 \vol Z
 =\sum_{\mathcal T\ {\rm tree}}\prod_j b_j^{\deg_{\mathcal T}(j)}
 =\Big(\prod_j b_j\Big)
    \sum_{(j_1,\ldots,j_{\ell-2})\in\{1,\ldots,\ell\}^{\ell-2}}
          \prod_{a=1}^{\ell-2}b_{j_a}
 =\Big(\prod_j b_j\Big)n^{\ell-2}.
\]
Indeed \(j\) occurs \(\deg_{\mathcal T}(j)-1\) times in the code. For \(\ell=1\) the last expression is one by \(b_1=n\). The preceding flag and diagonal calculations identify \(A_0\) on continuous regular tests. The sharp entry box contains the singular limit \(0\) for \(n\ge2\). Apply Proposition~\ref{localized} to its indicator, whose boundary has zero \(\eta_n\)-measure, and use Corollary~\ref{small} when some \(e_i>1\). By \eqref{ch},
\[
 C^\Box_{\Q,f}
 =\underbrace{\frac{\prod_i A_{p_i}^{e_i}}{L_n}}_{A_0}
  \underbrace{\frac{n^{\ell-2}\prod_i d_i^{e_i}}{\prod_i e_i!}}_{\vol(P)}
  \underbrace{J_{n,\R}}_{\eta_n(\Box)} .
\]
For \(n=1\) the same expression is one.
\end{proof}

\subsection{The residue of a rational simple order}
For monic irreducible \(p\in\Z[t]\) let \(F=\Q[t]/p\), \(R=\Z[t]/p\) and \(O=\Oa_F\). A fractional \(R\)-ideal \(I\) is a full \(\Z\)-lattice in \(F\) stable under \(R\). Put
\[
 I\sim I'\ \Longleftrightarrow\ I'=zI\ (z\in F^\times),\quad
 (R:I)=\{z\in F:zI\subset R\},\quad U_I=\{u\in F^\times:uI=I\}.
\]
Write \(v(I)\) for covolume in a fixed Lebesgue measure on \(F_\infty\), and \(A_p=\operatorname*{Res}_{s=1}\zeta_R(s)\).
\begin{prop}[Residue by lattice classes]\label{idealres}
For \(R=\Z[t]/p\) as above the residue \(A_p\) is
\begin{equation}\label{finiteres}
 A_p=\frac{\operatorname{Res}_{s=1}\zeta_F(s)}{h_F}
   \sum_{[I]} [O^\times:U_I]\frac{v(R)v(O)}{v(I)v((R:I))}.
\end{equation}
Here \(h_F\) is the class number of \(F\); the sum is over homothety classes of all full fractional \(R\)-ideals.
\end{prop}
\begin{proof}
Writing \(\theta=t\bmod p\), trace duality for a monogenic order gives \(R^\vee=p'(\theta)^{-1}R\) under the trace pairing. Multiplication by \(p'(\theta)\) therefore identifies the submodules counted by \(\zeta_R\) with integral \(R\)-ideals, preserving indices. For a fixed representative \(J\) of an \(O\)-ideal class and conductor \(\mathfrak f\), normalize \(OI=J\). Then
\[
 \mathfrak fJ=(\mathfrak f O)I\subset I\subset J,\qquad
 \#(J/\mathfrak fJ)<\infty,\qquad \operatorname{Stab}_{F^\times}(J)=O^\times .
\]
This gives a finite set of homothety classes.

In a fixed class define \(A_I(y)=\#\{J'\subset R:J'\sim I,\,[R:J']\le y\}\). Counting the parameters \(z\) gives
\[
 A_I(y)=\#\left\{z\in((R:I)\setminus0)/U_I:
        |Nz|\le y\,\frac{v(R)}{v(I)}\right\}
 \sim c_F[O^\times:U_I]\,\frac{v(R)}{v(I)v((R:I))}\,y .
\]
Here \(J'=zI\) and \([R:zI]=|Nz|v(I)/v(R)\). In the fundamental-region calculation of Lemma~\ref{shell}, passing from \(O^\times\) to \(U_I\) multiplies the region's volume by \([O^\times:U_I]\). Applying it to the lattice \((R:I)\) gives the asymptotic. For \(R=O\) each class contributes \(c_F/v(O)\), so
\[
 c_F=\frac{v(O)}{h_F}\operatorname{Res}_{s=1}\zeta_F(s),\qquad
 \lim_{s\downarrow1}(s-1)\zeta_R(s)
      =\lim_{y\to\infty}y^{-1}\sum_{[I]}A_I(y).
\]
Indeed if \(A(y)=\sum_{[I]}A_I(y)=ay+o(y)\), then for real \(s>1\)
\[
 (s-1)\zeta_R(s)=s(s-1)\int_1^\infty A(y)y^{-s-1}\,dy\longrightarrow a.
\]
For the error, split the integral at a fixed \(Y\), bound \(|A(y)/y-a|\) by its supremum on \(y\ge Y\), and let \(Y\to\infty\). Meromorphic continuation \cite{Yun} identifies the limit with the residue. Substituting \(c_F\) gives \eqref{finiteres}.
\end{proof}

\begin{example}[A nonmaximal order]
For \(p=t^2+4\) we have \(R=\Z[2i]\subset O=\Z[i]\). Normalize \(v(O)=1\). After rescaling a fractional ideal by a field element, its \(O\)-saturation is \(O\), since \(O\) has class number one. Then \(2O\subset I\subset O\), and \(I/2O\) generates \(O/2O\) as an ideal. It is either the whole space or a line spanned by a unit. The two unit lines are exchanged by \(i\). Thus \(I=O,R\) represent all classes up to homothety. For \(R\) the unit index is two and \((R:R)=R\); for \(O\) the unit index is one and \((R:O)=2O\). The respective summands in \eqref{finiteres}, in the order \(R,O\), are \(1,1/2\). Consequently \(A_p=3\pi/8\).

One can check the factor \(3/2\) directly at the conductor prime. An ideal in Hermite normal form in the basis \(1,2i\) has columns \((ck,0)^t,(cj,c)^t\), with \(c,k\ge1\), \(0\le j<k\), \(j^2+4\equiv0\bmod k\). Indeed stability of \((a,0)^t,(b,c)^t\) under \((x,y)\mapsto(-4y,x)\) gives \(c\mid a,b\) and then this congruence. The index is \(c^2k\). Write \(z=2^{-s}\). The numbers of roots modulo \(1,2,4,8\) are \(1,1,2,2\), and there is none modulo \(16\). As \(2\) ramifies in \(O\), this gives
\[
 \frac{\zeta_{R,2}(s)}{\zeta_{O,2}(s)}
  =\frac{1+z+2z^2+2z^3}{1+z}
  =1+2z^2,\qquad
 \left.\frac{\zeta_{R,2}(s)}{\zeta_{O,2}(s)}\right|_{s=1}=\frac32 .
\]
At odd primes \(R\) and \(O\) coincide locally.
\end{example}

\subsection{Why primes dividing resultants cancel}
For a monic squarefree \(a\in\Z[t]\), put
\[
 R_{a,l}=\Z_l[t]/a,\qquad
 \widetilde R_{a,l}=\text{its maximal order},\qquad
 \mathcal Z_l(a)=
       \left.\frac{\zeta_{R_{a,l}}(s)}{\zeta_{\widetilde R_{a,l}}(s)}\right|_{s=1}.
\]

\begin{lem}[Cancellation of the resultant in the local residue]\label{resultantcancel}
If \(a,b\in\Z[t]\) are monic, squarefree and coprime over \(\Q\), then for every finite prime \(l\),
\[
 \mathcal Z_l(ab)=\mathcal Z_l(a)\mathcal Z_l(b).
\]
Consequently, for distinct monic irreducible factors \(g=\prod_{i=1}^k p_i\),
\[
 \lim_{s\to1}(s-1)^k\zeta_{\Z[t]/g}(s)=\prod_{i=1}^k A_{p_i}.
\]
\end{lem}
\begin{proof}
Let \(\mathrm{OI}_l(a)\) be the integral orbital mass, with maximal compact centralizer of volume one. Yun's identity in the normalization of \cite[Theorem 10.1 and its proof]{CD} is
\[
 \mathcal Z_l(a)=\frac{\mathrm{OI}_l(a)}
                           {[\widetilde R_{a,l}:R_{a,l}]} .
\]
Fix stable intersection and projection lattices \(L_a,L_b\) in the two rational summands. Every lift of this pair has the form
\[
 L(u)=\{(x+uy,y):x\in L_a,\ y\in L_b\},\qquad
 u\in\Hom_{\Q_l}(L_b\otimes\Q_l,L_a\otimes\Q_l).
\]
Two parameters give the same lift exactly when their difference lies in
\(\Lambda=\Hom_{\Z_l}(L_b,L_a)\). If \(A,B\) are multiplication by \(t\), stability means \(Au-uB\in\Lambda\).
The Sylvester map \(S:u\mapsto Au-uB\) is invertible over \(\Q_l\) and preserves \(\Lambda\). Thus the set of lifts is
\[
 \mathcal L=S^{-1}\Lambda/\Lambda,\qquad
 |\mathcal L|=[\Lambda:S\Lambda]
       =|\det S|_l^{-1}=|\operatorname{Res}(a,b)|_l^{-1}.
\]
For the compact stabilizer \(U\) of \((L_a,L_b)\) and its stabilizer \(U_x\) of a lift,
\[
 \sum_{[x]\in U\backslash\mathcal L}\frac1{\vol U_x}
   =\frac1{\vol U}\sum_{[x]}[U:U_x]
   =\frac{|\operatorname{Res}(a,b)|_l^{-1}}{\vol U}.
\]
Summing over the base pairs gives
\(\mathrm{OI}_l(ab)=|\operatorname{Res}(a,b)|_l^{-1}\mathrm{OI}_l(a)\mathrm{OI}_l(b)\).
The exact sequence
\[
 0\longrightarrow R_{ab,l}\longrightarrow R_{a,l}\oplus R_{b,l}
 \xrightarrow{(u,v)\mapsto u-v}\Z_l[t]/(a,b)\longrightarrow0
\]
gives the same factor in the order index. Indeed the last module is the cokernel of multiplication by \(b\) on \(\Z_l[t]/a\), whose determinant is \(\operatorname{Res}(a,b)\). Since \(\widetilde R_{ab,l}=\widetilde R_{a,l}\oplus\widetilde R_{b,l}\),
\[
 [\widetilde R_{ab,l}:R_{ab,l}]
 =|\operatorname{Res}(a,b)|_l^{-1}
  [\widetilde R_{a,l}:R_{a,l}][\widetilde R_{b,l}:R_{b,l}].
\]
Divide the two identities and apply Yun's formula. Finally
\(\zeta_{\Z[t]/g}/\prod_i\zeta_{\Z[t]/p_i}\) is a finite product of local rational functions, and each has value one at \(s=1\). This proves the global residue identity.
\end{proof}

For example, let \(a=t\), \(b=t-p\), with \(p\) prime, and set \(x=p^{-s}\). In the basis \(1,t\), an ideal of \(\Z_p[t]/(t(t-p))\) has Hermite normal form
\[
 \left\langle(p^{h+k},0),(p^h j,p^h)\right\rangle,\qquad
 0\le j<p^k,\quad j(j+p)\equiv0\bmod p^k .
\]
The numbers of roots for \(k=0,1,2\) are \(1,1,p\). For every \(k\ge3\) the number is \(2p\). Thus
\[
 \frac{\zeta_{R_{t(t-p),p}}(s)}{\zeta_{\Z_p\times\Z_p}(s)}
  =\frac{1-x}{1+x}\left(1+x+px^2+\frac{2px^3}{1-x}\right)
  =1-x+px^2 .
\]
This equals one at \(s=1\) and exhibits the dependence on \(s\) away from the residue in Lemma~\ref{resultantcancel}.

\section{Split quadratic fibers and short eigenlines}\label{quadratic}
Suppose \(f=(t-\lambda)(t-\mu)\) with \(\lambda,\mu\in\Oa\). For a cyclic \(X\) the primary flag has one proper term after choosing the first eigenvalue. When \(\lambda\ne\mu\) either eigenline can occur first; for a nonscalar double root the eigenline is unique. An irreducible quadratic has no \(K\)-eigenline and is outside this argument. The scalar double-root matrix is counted separately. Put
\[
 a_K=\frac{2^{r_1}(2\pi)^{r_2}}{\sqrt{|\Delta_K|}},\qquad
 b_K=\frac{2^{r_1+r_2-1}h_K R_K}{w_K\zeta_K(2)}\,a_K^2 .
\]
Here \(r_1,r_2,h_K,R_K,w_K\) have their usual meanings (real places, complex places, class number, regulator and number of roots of unity), with \(R_K=1\) in unit rank zero.
\begin{thm}[Split quadratic fibers]\label{quadrank}
For \(\lambda,\mu\in\Oa\) fixed, \(m=\#\{\lambda,\mu\}\),
\[
 N^\Box_{(t-\lambda)(t-\mu)}(T)=m q a_K b_K\,T^q\log T+O_{K,\lambda,\mu}(T^q).
\]
\end{thm}

We count matrices on one labeled line, then subtract the overlap when two labeled lines are short.

\subsection{Counting on a fixed eigenline}
For \(X\in M_2(K)\), a \emph{labeled eigenline} is a pair \((P,\lambda)\) with \(P\in\mathbb P^1(K)\) and \((X-\lambda I)P=0\). For \(P=[u:v]\) define
\[
 I_P=u\Oa+v\Oa,\qquad M_w=\max(|u|_w,|v|_w),\qquad
 h(P)=\frac{\prod_{w\mid\infty}M_w^{q_w}}{N I_P}.
\]
Here \(NI_P\) is the absolute norm of the fractional ideal represented by \(u,v\).
Scaling \((u,v)\) by \(z\in K^\times\) multiplies the numerator and denominator by the same \(|N_{K/\Q}z|\). Thus \(h\) depends only on \(P\). Schanuel's estimate \cite{Schanuel} and partial summation give
\[
 \#\{P:h(P)\le x\}=b_Kx^2+O_K(x^{2-1/q}\log(2x)),\qquad
 \sum_{h(P)\le x}h(P)^{-2}=2b_K\log x+O_K(1).
\]
\begin{lem}[A labeled eigenline]\label{oneline}
For \(P\in\mathbb P^1(K)\), put \(h=h(P)\) and define
\[
\begin{aligned}
 M_{\lambda,\mu}(P,T)&=\#\{X\in M_2(\Oa):H_\Box(X)\le T,\
         X|_P=\lambda,\ X|_{K^2/P}=\mu\},\\
 E(P,T)&= M_{\lambda,\mu}(P,T)-a_KT^q/h^2 .
\end{aligned}
\]
Here \(X|_{K^2/P}\) denotes the induced endomorphism. Uniformly for \(T\ge2\), and for each fixed \(c\ge0\),
\[
 \begin{aligned}
 |E(P,T)|&\ll_{K,\lambda,\mu}
       1+(T h^{-2/q})^{q-1}+T^{q-1}h^{-2},\\
 \sum_{h(P)\le(T+c)^{q/2}} |E(P,T)|&\ll_{K,\lambda,\mu,c}T^q .
 \end{aligned}
\]
\end{lem}
\begin{proof}
Choose the representative \(I=I_P\) from a fixed list of ideal classes and balance the generators by units so that \(M_w\asymp_K h^{1/q}\). Define
\[
 B_P=\binom uv(-v,u),\qquad
 zB_P\in M_2(\Oa)\ \Longleftrightarrow\ zI^2\subseteq\Oa
                 \ \Longleftrightarrow\ z\in I^{-2}.
\]
These are exactly the integral maps killing \(P\) with image in \(P\). An integral splitting of the saturated line gives one lift \(X_0\); the full set is \(X_0+I^{-2}B_P\). At each place \(w\) let \(\Pi_{P,w}\) be orthogonal projection onto \(P_w\). The lift \(X_{*,w}=\mu_w I_2+(\lambda_w-\mu_w)\Pi_{P,w}\) satisfies \(\|X_{*,w}\|_{\max}\le c_0\), with \(c_0\) independent of \(P\). Since it has the two prescribed actions, \(X_{*,w}=X_{0,w}+\xi_w B_{P,w}\) for some \(\xi_w\). Thus
\[
 |z-\xi_w|M_w^2-c_0\le\|X_{0,w}+zB_{P,w}\|_{\max}
                          \le |z-\xi_w|M_w^2+c_0 .
\]
The local parameter region \(D_w(T)\) therefore satisfies, for \(T>c_0\),
\[
 \big\{|z-\xi_w|\le (T-c_0)/M_w^2\big\}\ \subseteq D_w(T)
 \ \subseteq\big\{|z-\xi_w|\le (T+c_0)/M_w^2\big\}.
\]
The centers are permitted to vary with \(P\). For the product at the middle radii the volume calculation is
\[
 \frac{2^{r_1}\pi^{r_2}T^q}{\prod_w M_w^{2q_w}}
 \cdot\frac{1}{2^{-r_2}|\Delta_K|^{1/2}(NI)^{-2}}
      =\frac{a_KT^q}{h^2}.
\]
A fixed-thickness boundary layer of a product of intervals and disks of radii \(\le L\) has volume \(O_K((1+L)^{q-1})\), also if one radius is shortened. A bounded fundamental cell for each of the finitely many lattices \(I^{-2}\) thus gives error \(O_K(1+(Th^{-2/q})^{q-1})\). Changing \(T\) to \(T\pm c_0\) adds \(O(T^{q-1}/h^2)\). The outer bound covers bounded \(T\) as well. On dyadic \(h\asymp H\), Schanuel's bound gives \(O_K(H^2)\) lines, so the total error is
\[
 \ll\sum_{\substack{H=2^j,\ j\ge0\\H\ll (T+c)^{q/2}}}
          \big(H^2+T^{q-1}H^{2/q}+T^{q-1}\big)\ll_{K,\lambda,\mu,c}T^q .
\]
\end{proof}

\begin{proof}[Proof of Theorem~\ref{quadrank}]
\stage{A double root}
Use \(B=(T+c)^{q/2}\), with a fixed \(c\) large enough for both eigenvalues. For a nonscalar matrix with double root, \(X-\lambda I_2=zB_P\) with \(z\ne0\); integrality gives \(|Nz|(NI)^2\ge1\), whence \(h(P)\le B\). There is a unique eigenline. Summing Lemma~\ref{oneline} over \(h\le B\) and omitting the scalar once per line costs \(O(T^q)\); restore the scalar once. The main term is
\[
 a_K T^q(2b_K\log B+O(1))=qa_Kb_K T^q\log T+O(T^q).
\]

\stage{Two distinct roots}
For different eigenvalues label their lines \(P,Q\). If \(p,p'\) span them with coordinate ideals \(I,J\), write \(\mathfrak d=(\det(p,p'))/(IJ)\). This is integral, and the matrix \(X-\mu I_2\) is \((\lambda-\mu)\) times the projection onto \(P\) along \(Q\). Its entry ideal is \((\lambda-\mu)/\mathfrak d\). Therefore \(1\le N\mathfrak d\le|N(\lambda-\mu)|\) and, using the maximum entry in the rank-one projection at each place,
\[
 \prod_{w\mid\infty}\|X-\mu I_2\|_{w,\max}^{q_w}
 = \frac{|N(\lambda-\mu)|}{N\mathfrak d}\,h(P)h(Q).
\]
It follows that \(h(P)h(Q)\le (T+c)^q\), so at least one labeled eigenline has height \(\le B\). Each choice of label contributes the same leading term just computed. We must subtract the overlap.

If both heights are at most \(B\) and \(h(P)=h\), the product on the left is \(O(hB)\). At one place the local norm is therefore \(O((T+c)(h/B)^{1/q})\). In Lemma~\ref{oneline} shorten the radius at that place by the factor \((h/B)^{1/q}\), up to a fixed constant. The uniformly bounded change of origin is absorbed since \(h\ge1\) and \((T+c)B^{-1/q}\ge1\). The volume contribution per line is \(O(T^q h^{-2}(h/B)^{1/q})\); the boundary errors still have total \(O(T^q)\). Sum over the finitely many places and over dyadic \(h\le B\). On \(2^{-j-1}B<h\le2^{-j}B\) there are \(O_K((2^{-j}B)^2)\) lines. Summing the shortened volumes and the boundary errors therefore gives
\[
 N_{\rm overlap}(T)\ll T^q\sum_{j\ge0}2^{-j/q}+O(T^q)
                 \ll_{K,\lambda,\mu}T^q .
\]
Inclusion--exclusion now yields the factor \(m=2\).
\end{proof}

For \(n=2\), the largest entry modulus of \(k E_{12} k^{-1}\) is \(\max(|k_{11}|^2,1-|k_{11}|^2)\). The real circle parameter and the complex Haar column give, respectively,
\[
 J_{2,\R}=\frac{8}{\pi}\int_0^{\pi/4}\sec^2\theta\,d\theta=\frac8\pi,\qquad
 J_{2,\C}=\pi\int_0^1\frac{dx}{\max(x,1-x)^2}=2\pi .
\]
For \(g=t(t-1)\), \(\Oa[t]/g=\Oa\times\Oa\), so \(\mathfrak S(g)=(\operatorname{Res}_{s=1}\zeta_K(s))^2\).
Inserting the angular values in \eqref{arith} for \(f=t^2\) and using
\(\operatorname{Res}_{s=1}\zeta_K(s)=a_K h_KR_K/w_K\) gives \(q a_K b_K\), as above. For \(\Q\) the two leading coefficients are \(24/\pi^2\) and \(48/\pi^2\); for \(\Q(i)\) they are \(\pi^3/(2\zeta_{\Q(i)}(2))\) and twice this value.
\subsection{The rational nilpotent cone}
\begin{prop}[An exact count and a secondary term]\label{cone}
Assume \(K=\Q\). For every integer \(H\ge0\),
\begin{equation}\label{totient}
 N_{t^2}^\Box(H)=1+8\sum_{m\le\sqrt H}\varphi(m)\lfloor H/m^2\rfloor .
\end{equation}
As \(H\to\infty\), the two-term expansion is
\begin{equation}\label{secondary}
 N_{t^2}^\Box(H)=\frac{24}{\pi^2}H\log H+
 \left(\frac{4(3\gamma-1)}{\zeta(2)}-\frac{8\zeta'(2)}{\zeta(2)^2}\right)H+O(H^{2/3}).
\end{equation}
Here \(\gamma\) is Euler's constant.
\end{prop}
\begin{proof}
Primitive unoriented pairs of maximum norm \(m\) give \(4\varphi(m)\) lines, with \(\varphi(1)=1\). On each line the nonzero integer parameter takes both signs and has modulus at most \(H/m^2\). Adding the zero matrix gives \eqref{totient}. The secondary term is the matrix version of Elkies's square-product estimate \cite{Elkies}. We give the calculation.
For the connection, a nilpotent matrix is \(\left(\begin{smallmatrix}x&b\\c&-x\end{smallmatrix}\right)\) with \(x^2=-bc\). Each ordered pair of positive integers \(\le H\) with square product gives four such matrices; the axes give \(4H+1\).

For \(A(x)=\sum_{b^2m\le x} b\), split at \(U=\lfloor x^{1/3}\rfloor\), \(V=\lfloor x/U^2\rfloor\). With \(L_m=\lfloor\sqrt{x/m}\rfloor\) the exact decomposition is
\[
 A(x)=\sum_{b\le U}b\lfloor x/b^2\rfloor+\tfrac12\sum_{m\le V}L_m(L_m+1)
       -\tfrac12 VU(U+1).
\]
Since \(U,V\asymp x^{1/3}\), summing \(b\lfloor x/b^2\rfloor=x/b+O(b)\) and \(L_m(L_m+1)=x/m+O(\sqrt{x/m})\) gives
\[
 \begin{aligned}
 A(x)&=x(\log U+\gamma)+\tfrac{x}{2}(\log V+\gamma)-\tfrac{x}{2}
      +O(x/U+U^2+x/V+\sqrt{xV}+UV)\\
 &=\tfrac{x}{2}\log x+\tfrac{3\gamma-1}{2}x+O(x^{2/3}).
 \end{aligned}
\]
Here \(U^2 V=x+O(U^2)\). Substituting \(\varphi(k)=\sum_{d\mid k}\mu(d)k/d\) into \eqref{totient} now gives
\[
 \frac{N_{t^2}^\Box(H)-1}{8}
   =\sum_{d\le\sqrt H}\mu(d)A(H/d^2)
   =\frac{H(\log H+3\gamma-1)}{2\zeta(2)}
      -\frac{H\zeta'(2)}{\zeta(2)^2}+O(H^{2/3}).
\]
Here \(\sum d^{-4/3}<\infty\). Extending the main sums to infinity costs \(O(\sqrt H\log(2H))\), and differentiating the reciprocal zeta series gives \(\sum_{d\ge1}\mu(d)\log d/d^2=\zeta'(2)/\zeta(2)^2\).
\end{proof}

\begin{cor}[The height of the kernel line]\label{directflag}
Let \(v_X\) be a primitive generator of the kernel of a nonzero integral nilpotent matrix \(X\), and set \(u_H(X)=\log\|v_X\|_2/\log H-1/2\). For each \(s\in\R\),
\[
 \lim_{H\to\infty}
 \frac{\#\{0\ne X\in M_2(\Z):X^2=0,\ H_\Box(X)\le H,\ u_H(X)\le s\}}
      {N_{t^2}^\Box(H)-1}
 =
 \begin{cases}
       0,&s\le-1/2,\\
       2s+1,&-1/2<s<0,\\
       1,&s\ge0.
 \end{cases}
\]
Thus the limit is the uniform law on \([-1/2,0]\).
\end{cor}
\begin{proof}
For fixed \(0<u\le1/2\), restrict the sum in \eqref{totient} to \(m\le H^u\). Removing the floors costs \(O(H^{2u})\). Writing \(m=dk\) in the divisor formula for \(\varphi\) gives
\[
 \sum_{m\le x}\frac{\varphi(m)}{m^2}
 =\sum_{d\le x}\frac{\mu(d)}{d^2}\sum_{k\le x/d}\frac1k
 =\frac{\log x}{\zeta(2)}+O(1),
\]
since \(\sum_{d\ge1}d^{-2}\log(2d)<\infty\).
The restricted count is \(8uH\log H/\zeta(2)+O(H)\). Its ratio to the total tends to \(2u\), the distribution function of the uniform law on \([0,1/2]\). Replacing the maximum coordinate by Euclidean length changes the logarithm by at most \(\log\sqrt2\). This also directly verifies the orientation of the chamber in Theorem~\ref{law}.
\end{proof}

\begin{prop}[The threshold for an integral translation]
For \(c\in\Z\) with \(a_c=|c|\ge1\),
\[
 N_{(t-c)^2}^\Box(H)=N_{t^2}^\Box(H)-4a_c
 \quad\bigl(H\ge (a_c-1)^2+1,\ H\in\Z\bigr).
\]
The least integer from which the equality holds at every height is \((a_c-1)^2+1\). The equality fails at \(H=(a_c-1)^2\) for every \(a_c\ge1\).
\end{prop}
\begin{proof}
For a primitive line \((u,v)\), with \(m=\max(|u|,|v|)\), direct inspection gives
\[
 H_\Box(cI_2+zB_P)=\max\{|z|m^2,\ |c-zuv|,\ |c+zuv|\}
                  =\max\{|z|m^2,\ a_c+|zuv|\}.
\] For \(|u|\ne|v|\), loss of a nonzero parameter implies
\[
 |z|m^2\le H<a_c+|z|m(m-1)
 \ \Longrightarrow\ |z|m\le a_c-1
 \ \Longrightarrow\ H\le (a_c-1)m\le(a_c-1)^2.
\]
Above this bound the scalar is present, and only the slopes \(1,-1\) lose points. Together they lose \(2\cdot2\#\{j\in\Z:H-a_c<j\le H\}=4a_c\) parameters. For sharpness take \(a_c\ge3\), \(m=a_c-1\), \(H=m^2\), and \((u,v,z)=(m,m-1,1)\). Then \(a_c+|uv|=m+1+m(m-1)=H+1\), so this third slope also loses a point. For \(a_c=1,2\) the preceding heights are \(0,1\); the trace excludes every translated matrix there, and \eqref{totient} gives \(1,9\).
\end{proof}

\appendix
\section{Series for the real and complex box factors}\label{angular}
For \(n=1\), \(J_{1,v}=1\). We give a convergent evaluation with an explicit error bound for \(n\ge2\).
Let \(d=q_vD_n,\ s=d/2,\ b=n^2\). Independently sample \(k\) from Haar probability on \(\mathcal K_v\) and \(U\) uniformly on the real Euclidean unit sphere of \(\mathfrak n_n(K_v)\). In row order set \(x_{ab}=|(kUk^{-1})_{ab}|^2\); denote these \(n^2\) coordinates by \(x_1,\ldots,x_b\). Then \(x_j\ge0\) and \(\sum_jx_j=1\). For \(m\ge1\) put
\[
 P_m=\sum_{|\alpha|=m}\binom m{\alpha}
        \left(\frac m{\max_j\alpha_j}\right)^{s}
        \mathbb E(x_1^{\alpha_1}\cdots x_b^{\alpha_b}).
\]
\begin{prop}[A convergent series for the box factor]
The moments in \(P_m\) are rational. With \(L=s b^{s+1}\),
\[
 J_{n,v}=\kappa_d\left(1+\sum_{h=1}^\infty(P_{4^h}-P_{4^{h-1}})\right),\qquad
 |J_{n,v}-\kappa_d P_m|\le\kappa_d L\sqrt{(1-1/b)/m}.
\]
The series converges absolutely.
\end{prop}
\begin{proof}
\stage{Error bound} Put \(F(x)=(\max_jx_j)^{-s}\) on the probability simplex. Conditional on \(x\), let \(M\sim\operatorname{Multinomial}(m,x)\). For a unit direction \(U\), the box cuts its radial line at \(r=(\max_j x_j)^{-1/2}\). With probability measure on directions the radial density is \(d\kappa_d\,r^{d-1}dr\), so
\[
 J_{n,v}=d\kappa_d\,\mathbb E\int_0^{(\max_jx_j)^{-1/2}}r^{d-1}\,dr
        =\kappa_d\mathbb EF(x),\qquad P_m=\mathbb EF(M/m).
\]
The second equality is the multinomial mass formula. On the simplex, \(\max_j x_j\ge1/b\). The mean value theorem gives
\[
 |F(x)-F(y)|\le sb^{s+1}|\max_j x_j-\max_j y_j|
        \le L\|x-y\|_2 .
\]
Therefore
\[
 \begin{aligned}
 \mathbb E(\|M/m-x\|_2^2\mid x)
   &=m^{-1}(1-\sum_j x_j^2)\le (1-1/b)/m,\\
 |J_{n,v}/\kappa_d-P_m|
   &\le L\,\mathbb E\|M/m-x\|_2
     \le L\sqrt{(1-1/b)/m}.
 \end{aligned}
\]
Since \(P_1=1\), the series telescopes, and
\(|P_{4^h}-P_{4^{h-1}}|\le 3L\sqrt{1-1/b}\,2^{-h}\) proves absolute convergence.

\stage{Rational moments} Condition on successive columns of \(k\). Let \(P\) project orthogonally onto the remaining \(a\)-dimensional subspace and let \(v\) be uniform on its unit sphere. In the complex case only balanced moments survive, and
\[
 \mathbb E\!\left(\prod_{j=1}^t v_{i_j}\overline{v}_{l_j}\right)=
 \frac{\sum_{\sigma\in S_t}\prod_{j=1}^t P_{i_j l_{\sigma(j)}}}{a(a+1)\cdots(a+t-1)} .
\]
In the real case odd moments vanish and
\[
 \mathbb E\prod_{j=1}^{2t}v_{i_j}
 =\frac{\displaystyle\sum_{\mathcal P}
           \prod_{\{u,w\}\in\mathcal P}P_{i_u i_w}}
         {\displaystyle\prod_{\ell=0}^{t-1}(a+2\ell)} ,
\]
where \(\mathcal P\) runs over pairings. Both moment identities follow by taking a projected standard Gaussian and factoring off its independent radius. At each step \(P=I-\sum_j k_jk_j^*\), the sum being over earlier columns. Integrate from the last column back to the first. Treat the real coordinates of the independent triangular \(U\) by the real sphere formula with \(P=I_d\). Expansion of the matrix products gives rational moments: the coefficients lie in \(\Q(i)\) and the final expectation is real. This supplies a finite recurrence at each \(m\).
\end{proof}

\section{An identity at finite height}\label{exact}
The affine dependence of the characteristic coefficients on one row or column was used by Habegger, Ostafe and Shparlinski \cite[Section 6.1]{HOS}. We apply that observation in coordinates on the last row.
Let \(T\ge0\), \(\mathcal D_T=\{u\in\Oa:\max_{v\mid\infty}|u|_v\le T\}\), and \(A\in\mathcal D_T^{(n-1)\times n}\). Write \(M_A(t)=t[I_{n-1}\ \ 0]-A\) and
\[
 c_j(t)=(-1)^{n+j}\det M_A(t)_{\widehat j}\quad(1\le j\le n),
\]
where \(\widehat j\) denotes deletion of column \(j\); an empty determinant is one. Choose an integral basis of \(\Oa\). For \(u(t)=\sum_{r=0}^{n-1}u_rt^r\), let \([u]\in\Z^{qn}\) be the concatenation of the integral coordinates of \(u_0,\ldots,u_{n-1}\).
The set \(\mathcal D_T\) is finite by discreteness of the Minkowski embedding.
\begin{prop}[Last-row coefficient extraction]
In \(qn\) Laurent variables, the entry-box count at every \(T\ge0\) is
\[
 N_f^\Box(T)=
 \sum_A [\,\mathbf z^{[t c_n-f]}\,]\;
       \prod_{j=1}^n\left(\sum_{u\in\mathcal D_T}\mathbf z^{[u c_j]}\right).
\]
\end{prop}
\begin{proof}
Laplace expansion along the last row yields
\[
 \det(tI-X)=\sum_j(t\delta_{jn}-u_j)c_j(t)
           =tc_n(t)-\sum_j u_jc_j(t).
\]
Since \(c_n\) is monic of degree \(n-1\),
\(\chi_X=f\ \Longleftrightarrow\ \sum_j[ u_j c_j]=[t c_n-f]\).
For each fixed \(A\) the product expands as
\[
 \sum_{(u_1,\ldots,u_n)\in\mathcal D_T^n}\mathbf z^{\sum_j[u_jc_j]} .
\]
Extracting the specified exponent counts its entire fiber with multiplicity one per tuple, including tuples with coincident exponents. Summing over \(A\) counts every matrix exactly once. For \(K=\Q\), \(\mathcal D_T=[-\lfloor T\rfloor,\lfloor T\rfloor]\cap\Z\).
\end{proof}

For example take \(K=\Q\), \(f=t^2(t-1)\). The noncyclic matrices are exactly the idempotents in the fiber. The trace eliminates the last diagonal entry. Checking the sum of principal two-by-two minors, the determinant and \(X^2=X\) then gives the following totals.
\[
\begin{array}{c|rrr}
T&\text{cyclic}&\text{noncyclic}&N_f^\Box(T)\\\hline
1&600&87&687\\
2&6864&447&7311\\
3&26256&999&27255
\end{array}
\]

For a complex-place check take \(K=\Q(i)\), \(n=2\), \(f=t^2\).
The same counts at \(T=1,2,3\) are \(25,105,233\). They can be reproduced by enumerating triples in the Gaussian disk \(\mathcal D_T\) satisfying \(x^2+bc=0\).
These tables check the finite-height identities. Numerical convergence to the repeated-primary asymptotic in degrees \(n\ge3\) has yet to be tested.

\section*{Statements and Declarations}
\noindent\textbf{Use of AI tools.}
The authors used LLMs for proof-strategy discussions, exploratory
computations, and editing. The authors take responsibility for the
arguments and the final text.

\end{document}